


  \newcount\fontset
  \fontset=1
  \def\dualfont#1#2#3{\font#1=\ifnum\fontset=1 #2\else#3\fi}

  \dualfont\bbfive{bbm5}{cmbx5}
  \dualfont\bbseven{bbm7}{cmbx7}
  \dualfont\bbten{bbm10}{cmbx10}
  \font \eightbf = cmbx8
  \font \eighti = cmmi8 \skewchar \eighti = '177
  \font \eightit = cmti8
  \font \eightrm = cmr8
  \font \eightsl = cmsl8
  \font \eightsy = cmsy8 \skewchar \eightsy = '60
  \font \eighttt = cmtt8 \hyphenchar\eighttt = -1
  \font \msbm = msbm10
  \font \sixbf = cmbx6
  \font \sixi = cmmi6 \skewchar \sixi = '177
  \font \sixrm = cmr6
  \font \sixsy = cmsy6 \skewchar \sixsy = '60
  \font \tensc = cmcsc10
  
  \font \titlefont = cmbx10 scaled \magstep1
  \scriptfont \bffam = \bbseven
  \scriptscriptfont \bffam = \bbfive
  \textfont \bffam = \bbten


  \newskip \ttglue

  \def \eightpoint {\def \rm {\fam0 \eightrm }%
  \textfont0 = \eightrm
  \scriptfont0 = \sixrm \scriptscriptfont0 = \fiverm
  \textfont1 = \eighti
  \scriptfont1 = \sixi \scriptscriptfont1 = \fivei
  \textfont2 = \eightsy
  \scriptfont2 = \sixsy \scriptscriptfont2 = \fivesy
  \textfont3 = \tenex
  \scriptfont3 = \tenex \scriptscriptfont3 = \tenex
  \def \it {\fam \itfam \eightit }%
  \textfont \itfam = \eightit
  \def \sl {\fam \slfam \eightsl }%
  \textfont \slfam = \eightsl
  \def \bf {\fam \bffam \eightbf }%
  \textfont \bffam = \eightbf
  \scriptfont \bffam = \sixbf
  \scriptscriptfont \bffam = \fivebf
  \def \tt {\fam \ttfam \eighttt }%
  \textfont \ttfam = \eighttt
  \tt \ttglue = .5em plus.25em minus.15em
  \normalbaselineskip = 9pt
  \def \MF {{\manual opqr}\-{\manual stuq}}%
  \let \sc = \sixrm
  \let \big = \eightbig
  \setbox \strutbox = \hbox {\vrule height7pt depth2pt width0pt}%
  \normalbaselines \rm }


 \def \Headlines #1#2{\nopagenumbers
    \advance \voffset by 2\baselineskip
    \advance \vsize by -\voffset
    \headline {\ifnum \pageno = 1 \hfil
    \else \ifodd \pageno \tensc \hfil \lcase {#1} \hfil \folio
    \else \tensc \folio \hfil \lcase {#2} \hfil
    \fi \fi }}

  \def \Date #1 {\footnote {}{\eightit Date: #1.}}


  \def \lcase #1{\edef \auxvar {\lowercase {#1}}\auxvar }

  \def \goodbreak {\vskip0pt plus.1\vsize \penalty -250 \vskip0pt
plus-.1\vsize }

  \newcount \secno  \secno = 0
  \newcount \stno   \stno =0
  \newcount \eqcntr \eqcntr=0

  \def \seqnumbering {\global \advance \stno by 1 \global \eqcntr=0
    \number \secno .\number \stno }

  \def\eqmark#1{\global \advance\eqcntr by 1 
    \edef\a{\number \secno .\number\stno.\number\eqcntr}
    \eqno {(\a)}
    \global \edef #1{\a}\track{showlabel}{*}{#1}}

  \def\section #1{\global\def\SectionName{#1}\stno = 0 \global
\advance \secno by 1 \bigskip \bigskip \goodbreak \centerline {\tensc 
\number \secno .\enspace #1.}\bigskip \noindent \ignorespaces}

  \long \def \sysstate #1#2#3{\medbreak \noindent {\bf \seqnumbering
.\enspace #1.\enspace }{#2#3\vskip 0pt}\medbreak }
  \def \state #1 #2\par {\sysstate {#1}{\sl }{#2}}
  \def \definition #1\par {\sysstate {Definition}{\rm }{#1}}
  \def \remark #1\par {\sysstate {Remark}{\rm }{#1}}


  \def \proof {\medbreak \noindent {\it Proof.\enspace }}
  \def \proofend {\ifmmode \eqno \square \else \hfill \square
\looseness = -1 \medbreak \fi }

  \def \$#1{#1 $$$$ #1}
  \def\=#1{\buildrel #1 \over =}

  \def\iItem {\smallskip}
  \def\Item #1{\smallskip \item {#1}}
  \newcount \zitemno \zitemno = 0
  \def\izitem {\zitemno = 0}
  \def\zitem {\global \advance \zitemno by 1 \Item {{\rm(\romannumeral
\zitemno)}}}

  \newcount \footno \footno = 1
  \newcount \halffootno \footno = 1
  \def\footcntr {\global \advance \footno by 1
  \halffootno =\footno
  \divide \halffootno by 2
  $^{\number\halffootno}$}
  \def\fn#1{\footnote{\footcntr}{\eightpoint#1}}


  \def \N {{\bf N}}
  \def \C {{\bf C}}
  \def \<{\left \langle }   \def \<{\langle }
  \def \>{\right \rangle }  \def \>{\rangle }
  
  \def \ds{\displaystyle}
  \def \and {\hbox {,\quad and \quad }}
  \def \calcat #1{\,{\vrule height8pt depth4pt}_{\,#1}}
  
  \def \for #1{,\quad \forall\,#1}
  \def \square {\hbox {$\sqcap \!\!\!\!\sqcup $}}
  \def \crossproduct {{\hbox {\msbm o}}}
  \def \stress #1{{\it #1}\/}
  \def \inv {^{-1}}
  \def \* {\otimes}



  \def \ifn #1{\expandafter \ifx \csname #1\endcsname \relax }
  \def \track #1#2#3{\ifn{#1}\else {\tt [#2\string #3] }\fi}
  \def \cite #1{{\rm [\bf #1\track{showref}{\#}{#1}\rm]}}
  \def \scite #1#2{{\rm [\bf #1\track{showref}{\#}{#1}{\rm \hskip 0.7pt:\hskip 2pt #2}\rm]}}
  \def \label #1{\global \edef #1{\number \secno \ifnum \number \stno
                 = 0\else .\number \stno \fi }\track{showlabel}{*}{#1}}
  \def \lcite #1{(#1\track{showcit}{$\bullet$}{#1})}

  \newcount \bibno \bibno =0

  \def \newbib #1#2{\global\advance\bibno by 1 \edef#1{\number\bibno}}
  \def \newbib #1#2{\edef #1{#2}}

  \def \bibitem #1#2#3#4{\smallskip \item {[#1]} #2, ``#3'', #4.}
  \def \references {
    \begingroup 
    \bigskip \bigskip \goodbreak 
    \eightpoint 
    \centerline {\tensc References}
    \nobreak \medskip \frenchspacing }



  \def\lq#1#2{#1/#2}
  \def\rq#1#2{#2\backslash#1}
  \def\Rq#1#2{\displaystyle{#1\over #2}}

  \def\dq#1#2#3{#2\backslash#1/#3}

  \def\Z{{\bf Z}} 
  \def\R{{\bf R}}
  \def\H#1{H^{#1}}
  
  \def\HH#1{H\cap \H{#1}}

  \def\conj#1{\overline{#1}}

  \def\Rat{{\bf Q}}
  \def\P{{\cal P}}
  \def\AP{{\bf A}_{\P}}
  \def\HP{H_{\P}}

  \def\AXB#1#2{\pmatrix{1 & #1 \cr 0 & #2}}
  \def\axb#1#2{[\,#1,#2\,]}
  \def\axb#1#2{\AXB{#1}{#2}}
  \def\mod{{\rm mod}}
  \def\Ker{{\rm Ker}}

  \def\F#1{\Fd\left(#1\right)}
  \def\Fd{F}
  
  \def\Q#1{Q_{#1}}   \def\Q#1{Q_{_{#1}}}
  \def\d#1{\delta_{#1}}
  \def\dH#1{\d{H#1}}
  \def\dpH#1{\delta'_{H#1}}
  
  \def\ip[#1,#2]{\<\kern-2pt\<\,#1,#2\,\>\kern-2pt\>}
  \def\letterfors{\sigma}
  \def\letterfort{\tau}
  \def\normal{\mathrel{\raise 1.3pt \hbox{$\triangleleft$} \kern -4.5 pt
\vrule height 0pt depth 0.3pt width 3.9pt \kern 2pt}}

  \def\p#1{P\kern-2.5pt_{_{#1}}}
  \def\wt{\widetilde}
  \def\pit{\widetilde\pi}
  \def\q#1{q_{#1}}      
  \def\qalt#1{q^2_{#1}} 

  \def\s#1{\letterfors_{#1}} \def\calcat#1{(#1)} \def\e#1{e_{#1}} 
  \def\t#1{\letterfort_{#1}}    
  \def\Calcat#1{|_{#1}}
  \def\letterforr{\xi}
  \def\r#1{\letterforr(#1)}
  \def\si#1{\s{#1\inv}}
  \def\ei#1{\e{#1\inv}}
  \def\w#1#2{w_{#1,#2}}

\def\Hecke{{\cal H}}
\def\HGH{\Hecke(G,H)}
\def\B{{\cal B}}

  \def \prelabel#1#2{\global\edef #1{#2}\lcite{#2}\track{showlabel}{*}{#1}}
  \def \relabel #1{\label\a\track {showlabel} {*}{#1}\ifx#1\a\else
\message{*** RELABEL CHECKED FALSE ***} RELABEL CHECKED FALSE,
EXITING. \end \fi}


  \newbib{\ALR}{ALR}
  \newbib{\BostConnes}{BC}
  \newbib{\Brenken}{B}
  \newbib{\cross}{DE1}
  \newbib{\ExelMisha}{DE2}
  \newbib{\partial}{DEP}
  \newbib{\newpim}{E1}
  \newbib{\TPA}{E2}
  \newbib{\Inverse}{E3}
  \newbib{\ortho}{E4}
  \newbib{\infinoa}{EL1}
  \newbib{\kms}{EL2}
  \newbib{\topfree}{ELQ}
  \newbib{\Willis}{GW}
  \newbib{\Hall}{H}
  \newbib{\QuiggKaliLand}{KLQ}
  \newbib{\Kassel}{Ka}
  \newbib{\Krieg}{Kr}
  \newbib{\Laca}{L}
  \newbib{\LacaFrank}{LF}
  \newbib{\LaLaOne}{LL1}
  \newbib{\LaLaTwo}{LL2}
  \newbib{\LROne}{LR1}
  \newbib{\LRTwo}{LR2}
  \newbib{\Ward}{W}

  \def\titletext{HECKE ALGEBRAS FOR PROTONORMAL SUBGROUPS}

  \Headlines
  {\titletext}
  {R.~Exel}

  \null
  \centerline{\titlefont \titletext}
  \bigskip
  \centerline{{\tensc Ruy Exel}\footnote{$^\star$}{\eightpoint
Partially supported by CNPq.}}
  \bigskip

  \Date{25 May 2005}

  \midinsert 
  \narrower \narrower \narrower 
  \eightpoint \noindent {\tensc Abstract}.  We introduce the term
\stress{protonormal} to refer to a subgroup $H$ of a group $G$ such
that for every $x$ in $G$ the subgroups $x\inv Hx$ and $H$ commute as
sets.  If moreover $(G,H)$ is a Hecke pair we show that the Hecke
algebra $\HGH$ is generated by the range of a canonical partial
representation of $G$ vanishing on $H$.  As a consequence we show that
there exists a maximum C*-norm on $\HGH$, generalizing previous
results by
  Brenken,
  Hall,
  Laca, Larsen,
  Kaliszewski, Landstad and Quigg.
  When there exists a normal subgroup $N$ of $G$, containing $H$ as a
normal subgroup, we prove a new formula for the product of the
generators and give a very clean description of $\HGH$ in
terms of generators and relations.  We also give a description of
$\HGH$ as a crossed product relative to a twisted partial action of
the group $\lq GN$ on the group algebra of $\lq NH$.
  Based on our presentation of $\HGH$ in terms of generators and
relations we propose a generalized construction for Hecke algebras in
case $(G,H)$ does not satisfy the Hecke condition.
  \endinsert

  \section{Introduction}
  After the pioneering work of Bost and Connes \cite{\BostConnes},
several authors started a systematic investigation of C*-algebras
obtained as completions of Hecke algebras.  It was quickly realized
  \cite{\ALR},
  \cite{\Brenken},
  \cite{\LROne}
  that the Hecke C*-algebra which plays the central role in
\cite{\BostConnes} may be successfully
described as the crossed product algebra  relative to  a  semigroup of endomorphisms,
prompting a large interest in the application of crossed product
techniques to study Hecke algebras.  See also
  \cite{\LaLaOne},
  \cite{\LaLaTwo} and
  \cite{\LacaFrank}.

  The objective of the present paper is to study Hecke algebras from a
similar point of view, namely the theory of partial group
representations \scite{\Inverse}{6.2} and twisted partial crossed
products \cite{\TPA}, \cite{\ExelMisha}. See also
  \cite{\cross},
  \cite{\partial},
  \cite{\newpim},
  \cite{\ortho},
  \cite{\infinoa},
  \cite{\kms}, and
  \cite{\topfree}.

If $H$ is a subgroup of a group $G$ recall that $(G,H)$ is said to be
a \stress{Hecke pair} if for each $x$ in $G$ the double coset $HxH$ is
the disjoint union of \stress{finitely many} right cosets; the number
of right cosets involved usually being denoted in the literature by
$R(x)$.  Some authors \cite{\BostConnes} also express the fact that
$(G,H)$ is a Hecke pair by saying that $H$ is an
\stress{almost-normal} subgroup of $G$.

Given a Hecke pair and a field $F$ one defines the Hecke algebra
$\HGH$ as being the $F$-algebra formed by all $F$-valued
finitely supported functions on the double coset space $\dq GHH$,
under a certain convolution product.

This algebra is therefore obviously linearly generated by the simplest
possible functions $1_{HxH}$ (the characteristic function of the
singleton $\{HxH\}$), where $x$ ranges in a family of representatives
for the double coset space $\dq GHH$.  For technical purposes we assume  that
the characteristic of $F$ is zero and use
  $$
  \s x = {1\over R(x)}\ 1_{HxH}
  \for x\in G.
  $$
  The starting point for our research can be subsumed by the question
as to what extent the map 
  $$
  x\in G \mapsto \s x \in \HGH
  $$
  is a group representation.  The most naive form of this question,
namely expecting that $\letterfors$ be a genuine group representation,
is not too interesting since this holds if and only if $H$ is a normal
subgroup of $G$, in which case $\HGH$ trivializes, being just the
group algebra of the quotient group.
  
This is where the theory of partial group representations comes into
play. Recall that a partial representation of a group $G$ in a unital
algebra $A$ is a map
  $
  u: G\to A,
  $
  such that $u(1)=1$, and the usual group law ``$u(xy)=u(x)u(y)$''
holds after it is left-multiplied by $u(x\inv)$ or right-multiplied by
$u(y\inv)$.  See \prelabel{\DefPreps}{2.1} below for a detailed
definition.
  
It therefore makes sense to ask when is $\letterfors$ a partial
representation.  Unfortunately the answer is again negative for many
Hecke pairs, including most examples associated to the modular group
${\sl SL}_2(\Z)$ discussed e.g. in \cite{\Krieg}.

But, on the fortunate side, there are interesting examples for which the
answer is affirmative.  Among these is the Hecke pair appearing in the
already mentioned work by Bost and Connes \cite{\BostConnes}, as well
as some, but not all, Hecke pairs appearing in the papers that came in
its wake.

Our first major effort is therefore directed at classifying the Hecke
pairs for which $\letterfors$ is a partial representation.  In pursuit
of this goal I have been led to considering a very weak normality
property: let us say that a subgroup $H$ of a group $G$ is
\stress{protonormal} if for every $x$ in $G$ the conjugate subgroup
  $$ 
  \H x=x\inv H x
  $$
  commutes with $H$ in the sense that the products of sets $\H x H$
and $H\H x$ coincide. 

There is not much in the literature about this property except for
some conditions for subnormality based on it for finite groups;  see
\cite{\Ward} and the references given there for more details. Also,
it seems to me that
this condition is related to Drinfeld's notion of \stress{quantum
double} (see \scite{\Kassel}{Chapter IX}) and perhaps it is interesting to
explore this relationship further, a task I have not undertaken.

In what I believe is the main contribution of the present work,
Theorems \prelabel{\MainPrep}{8.1} and
\prelabel{\NecessaryCondForPrep}{8.2} prove that $\letterfors$ is a
partial representation if and only if $H$ is protonormal.

  It is elementary to check, for instance, that for the Hecke pair in
\cite{\BostConnes} this condition is fulfilled.
  That Hecke pair is in fact a ``bit more
normal than protonormal''.  Recall from \cite{\Ward} that the subgroup
$H\subseteq G$ is said to be $n$-subnormal if there exists a normal
chain
  $$
  H=H_0\normal H_1\normal H_2\normal\ldots \normal H_n=G,
  $$
  of length $n$.
  If $H$ is 2-subnormal in $G$ then for every $x$ in $G$ and every $h$
in $H$ one has that $x\inv h x\in H_1$, so that $\H x\subseteq H_1$.
Since $H$ is normal in $H_1$ one has that $yH=Hy$, for all $y\in \H
x$, and hence $\H x H=H\H x$.  In other words, 2-subnormal subgroups
are necessarily protonormal.

Given the relevance of 2-subnormal subgroups in this work we shall
call these simply \stress{subnormal}.

The first proof I found of the fact that $\letterfors$ is a partial
representation assumed that $H$ is subnormal, but in trying to prove
that subnormality is a necessary condition for $\letterfors$ being a
partial representation I could only prove that $H$ must be
protonormal.  So the desire to generalize to protonormal groups came
naturally.  Having been born in such a roundabout way, I wonder how
relevant the notion of protonormal subgroups will ever be.
  After fiddling a bit with this notion I was able to find a
curious example of a Hecke pair $(G,H)$ such that $H$ is protonormal
in $G$ but not subnormal.  This seems to be based on the exceptional
properties of the prime number 2.  The reader will find the relevant
results in \prelabel{\ExampleNotSubnormal}{14.2} and
\prelabel{\ExampleIsProtonormal}{14.3} below.

Although Brenken does not mention the word ``subnormal'' in
\cite{\Brenken}, he often works under the assumption that there exists
a normal subgroup $N$ of $G$, containing $H$, and contained in the
normalizer of $H$.  Clearly the existence of such a subgroup $N$ is
tantamount to the fact that $H$ is subnormal in $G$.  Our results
therefore generalize some of the results in \cite{\Brenken}.  See also
\scite{\QuiggKaliLand}{Theorem 8.5}.

When the base field $F$ is equipped with an involution (as defined
precisely in the next section) such as the usual involution on the
field of complex numbers, Hecke algebras over $F$ can be made into
*-algebras by considering the involution (as in \cite{\BostConnes})
  $$
  f^\#(x) = \conj {f(x\inv)}
  \for x\in G,
  $$
  for all finitely supported functions $f$  on  $\dq GHH$.

  As in the theory of unitary group representations, most
partial representations of interest taking values in a  *-algebra satisfy
the identity
  $$
  u(x)^*=u(x\inv).
  $$
  Since $\HGH$ is a *-algebra it is natural to ask if this is the case
for $\letterfors$.  The answer is no but there exists another
involution on $\HGH$ with respect to which $\letterfors$ satisfies the
condition above.  This involution was already used in
\cite{\QuiggKaliLand} and is defined by
  $$
  f^*(x) = \Delta(t\inv)\ \conj{f(x\inv)},
  $$
  where $\Delta(x) = R(x)/R(x\inv)$ (recall from \scite{\Krieg}{I.3.6}
that $\Delta$ is a group homomorphism).  En passant, the similarity
with the formula for the adjoint in the C*-algebra of a locally
compact group given in \scite{Pedersen}{7.1}, $\Delta$ representing
the modular function there, is nicely explained in
\scite{\QuiggKaliLand}{Section 4}.

If our field allows for taking square roots, or more precisely if there
exists a multiplicative map $\lambda$ from $G$ to $F$ such that
$\lambda(x)^2=\Delta(x)$ for all $x$ in $G$ (which is clearly the case
if $F$ is the field of complex numbers) then the two involutions are
isomorphic (see \prelabel{\CompareStars}{5.6}).  Assuming that $F={\bf
C}$ and that $H$ is protonormal observe that $\letterfors$ being a partial
representation gives
  $$
  \s x \s x^* \s x = 
  \s x \si x \s x =
  \s x \s{x\inv x} =
  \s x \s 1 =
  \s x,
  $$
  so that any *-representation of $\HGH$ on a Hilbert space must send
the generating elements $\s x$ to partial isometries, and hence to
operators with norm no bigger than 1.  Therefore, for every $a\in\HGH$
the supremum of $\|\pi(a)\|$, as $\pi$ range in the collection of all
*-representations of $\HGH$, is a finite real number.  This supremum
defines a C*-norm on $\HGH$ which is obviously the maximum among all
such.  This solves a problem which has been addressed by many authors
  \scite{\Brenken}{Proposition 2.8},
  \scite{\Hall}{Corollary 4.6},
  \scite{\LaLaOne}{Proposition 1.4},
  \scite{\QuiggKaliLand}{Theorem 8.5}.

Our next main effort has got to do with the formula for the product
$\s x\s y$.  Since $\HGH$ is linearly generated by the $\s x$, its
multiplication operation is completely described be the ``structure
constants'' $\lambda_{x,y}^z$ implicitly defined by
  $$
  \s x \s y = \sum_{HzH\in\dq GHH} \lambda_{x,y}^z \s z.
  $$
  The reader will find formulas for these constants in
\scite{\Krieg}{I.4.4} and \cite{\QuiggKaliLand}.

Based on the techniques we developed we were able to find a
significant simplification for these formulas under the hypothesis
that $H$ is subnormal.  In fact, given $x$ and $y$ in $G$ it is easy
to show, based on the defining property of Hecke pairs, 
that $HxHyH$ is the disjoint  union of finitely many double cosets,
say
  $$
  HxHyH = \dot{\bigcup_{1\leq i\leq n}} \ Hz_iH.
  $$
  We prove in Theorem
\prelabel{\MainTheoremOnHeckeProduct}{10.2} that
  $$
  \s x\s y =
  {1\over n}\  \sum_{i=1}^n \s{z_i}.
  \eqno{(\dagger)}
  $$
  Thus, viewing the Hecke algebra as the algebra generated by double
cosets, as some authors have it, we see that the product of the double
cosets $HxH$ and $HyH$ in the Hecke algebra is very closely related to
the set theoretic product $HxH{\cdot} HyH$ in $G$: the former is precisely the
average of the double cosets contained in the latter.  In particular
there is no mention to right or left cosets as in most other product formulas.
  
Based on concrete examples we were able to determine that $(\dagger)$
does not hold in general.  It is therefore an interesting question
(see \prelabel{\OpenQuestion}{10.3}) to precisely determine for which
Hecke pairs does this hold.  I would very much like to know, for
instance, whether $(\dagger)$ holds for protonormal subgroups, a
question I have tried to solve without success.

  Back to the subnormal situation a straightforward but interesting
aspect about the above product formulas is that they emcompass the
whole algebraic structure of $\HGH$.  Precisely speaking we show in
Theorem \prelabel\Relations{10.5} that $\HGH$ is the universal
$F$-algebra generated by symbols $\{\s x\}_{x\in G}$ under relations
$(\dagger)$.  This should be compared to other descriptions of Hecke
algebras in terms of generators and relations, e.g.
\scite{\BostConnes}{Proposition 18}, \scite{\Brenken}{Theorem 3.10}
and \scite{\LaLaOne}{Theorem 1.9}.

Motivated by \cite{\LROne} we then take up the problem of describing
Hecke algebras as crossed products.  In order to describe our results
in that direction let $(G,H)$ be a Hecke pair and suppose that there
exists a subgroup $N$ of $G$ such that
  $
  H\normal N\normal G.
  $
  Clearly this implies that $H$ is subnormal in $G$.  

One may motivate the desire to describe $\HGH$ as a crossed product as
follows: since this algebra arises as an attempt to make sense of the
group algebra of the quotient $G/H$ (which is only a group if $H$ is
normal in $G$), it should be obtained somehow as a product of $G/N$ by
$N/H$.

In \cite{\LROne} and \cite{\Brenken} it is assumed that $G$ is a
semidirect product $N\crossproduct K$ for some group $K$ (in which case
$K$ is clearly isomorphic to $G/N$) and it is proved, under suitable
hypothesis, that $\HGH$ is a crossed product of the group algebra of
$H/N$ by a semigroup of endomorphisms somehow based on $K$.

Our description of $\HGH$ as a crossed product is based not on the
theory of crossed products by endomorphisms, but on the recent theory
of crossed products by partial actions \cite{TPA}, \cite{\ExelMisha}
briefly described in the next section.
  Precisely because this theory allows for a ``twisting cocycle'' we
do not need to assume a semidirect product structure on $G$.  Our main
result in that direction, Theorem
\prelabel{\MainTheoremOnCrossedProduct}{11.9}, then provides an
isomorphism
  $$
  \HGH \simeq  \F{\Rq NH}\crossproduct \Rq GN 
  $$
  where $\F{\lq NH}$ is the group algebra of the quotient group $\lq
NH$ and the crossed product is with respect to a certain twisted
partial action of the quotient group $\lq GN$ on $\F{\lq NH}$.  
  
If $G$ does have a semidirect product structure we may get rid of the
cocycle, a result we prove in Corollary in
\prelabel{\UntwistedCrossProduct}{11.10}.

  We should mention that \cite{\LaLaOne} proves a similar result in
which $H\subseteq N\normal G$, but $H$ is not supposed to be normal in
$N$ (the Hecke algebra for the pair $(N,H)$ replaces the group algebra
$\F{\lq NH}$), although it is still assumed that $G$ is a semidirect
product.  A common generalization therefore seems a worthwhile
project.

  Another interesting crossed product description for Hecke algebras,
based on Green's twisted crossed products, may be found in
\cite{\QuiggKaliLand}.

Perhaps an advantage of the partial crossed product description over
endomorphism crossed products is that we need not care at all about
the existence of certain generating subsemigroups required in 
\scite{\LaLaOne}{Theorem 1.9} or \scite{\Brenken}{Theorem 3.12}.

Recall that our description of Hecke algebras in terms of generators
and relations in \lcite{\Relations} refers to the decomposition of $HxHyH$
as a disjoint union of finitely many double cosets.
  One could then be tempted to do away with the Hecke condition, namely
that every double cosets contains finitely many right cosets, and
introduce a generalized condition by saying that  $(G,H)$ is a
\stress{pseudo Hecke pair} if for every $x$ and $y$ in $G$ one has
that $HxHyH$ is made out of  finitely many double cosets.
Unfortunately though,  at least in the case of a subnormal $H\subseteq G$, one may
prove with the aid of Propositions \prelabel{\DoubleInDoubleProp}{10.1}
and \prelabel{\DoubleInRightProp}{3.2} that every pseudo Hecke
pair is a true Hecke pair and vice-versa, so no extension of the usual
concept is obtained.

Nevertheless, based on some insight provided by Cuntz-Krieger algebras
for infinite matrices \cite{\infinoa}, we risk to introduce a
generalized Hecke algebra for a group-subgroup pair $(G,H)$ which does
not satisfy the Hecke condition. See Definition
\prelabel\GeneralizedHeckeAlgebra{13.1}.  Not having taken a single
step in the description of the beast thus brought into existence, we
at least give an example which might be of interest to some.

I would like to express my gratitude to a number of colleagues who, in
a way or another, knowingly or not, were instrumental for the
completion of this work.  Those include, but are not limited to,
M.~Dokuchaev, D.~Evans, and A.~Zalesski who, over a short lunch,
showed me a smooth path to the basic theory of Hecke algebras.

\goodbreak \noindent For the readers's convenience this work is divided
up into the following sections:

\newcount\indxcount\indxcount=0
\def\indx#1{\par \global\advance\indxcount by 1
\itemitem{\number\indxcount.} {\sl #1.}}

  \medskip
  \indx{Introduction}
  \indx{Generalities about partial representations}
  \indx{Generalities about Hecke pairs}
  \indx{The Hecke algebra}
  \indx{$*$-algebra structure} 
  \indx{Commuting subgroups}
  \indx{Protonormal subgroups}
  \indx{The canonical partial representation}  
  \indx{Generalities about subnormal groups}
  \indx{A formula for the product and relations for the Hecke algebra}
  \indx{Hecke algebra as a crossed product}
  \indx{Hecke C*-algebras}
  \indx{A possible generalization of Hecke algebras}
  \indx{An example}

\section{Generalities about partial representations}
  Let $F$ be a field of characteristic zero\fn{One may perhaps
generalize our results to other fields by tracking that its
characteristic does not divide the order of certain coset spaces to be
considered later.}. We will assume that $F$ has a
\stress{conjugation}, that is, an involutive automorphism
  $$
  z\in F \mapsto \bar z \in F,
  $$
  which will be fixed form now on.  In the absence of a more interesting
conjugation one could take the identity map by default.  Clearly when
$F$ is the field of complex numbers the conjugation of choice should
be the standard one.

A map $\phi:U\to V$ between $F$-vector spaces $U$ and $V$ will be
called  \stress{conjugate-linear}  when it is additive and
$\phi(\lambda u)=\conj \lambda \phi(u)$ for all $\lambda\in F$ and
$u\in U$.

A \stress{*-algebra} is by definition an algebra $A$ over $F$ equipped
with an involution
  $$
  a\in A \mapsto a^*\in A
  $$
  which is conjugate-linear and such that $(ab)^*=b^*a^*$, for all $a$
and $b$ in $A$.

Whenever we speak of the group algebra $\F G$, for a given group $G$, we
will think of it as a *-algebra with the the unique involution such
that
  $$
  (\d t)^* = \d{t\inv}
  \for t\in G,
  $$
  where $\d t$ refers to the group element $t$ interpreted as an
element of $\F G$.

By a \stress{sesqui-linear} form on an $F$-vector space
$V$ we will  mean a function $$\phi:V\times V \to F,$$ which is
linear in the first variable and conjugate-linear in the second
variable.  We will say that $\phi$ is a \stress{hermitian} form if
$\phi$ moreover satisfies
  $$
  \phi(u,v)=\conj{\phi(v,u)}
  \for u,v\in V.
  $$
  A non-degenerate hermitian form will be one for which
  $$ 
  \big(\forall v\ \phi(u,v)=0) \Rightarrow u=0.
  $$

  We shall now list a few definitions of relevance to the later
sections for the convenience of the reader.  See the references given
for more information.

\definition
  \relabel\DefPreps
  \cite{\Inverse, \partial}.\quad
  A partial representation of a group $G$ in a unital algebra $A$ is,
by definition, a map $\letterfors:G\to A$ such that
  \izitem
  \zitem $\s 1=1$,
  \zitem $\si x \s x \s y = \si x \s{xy}$,
  \zitem $\s x \s y \si y = \s{xy} \si y$,
  \medskip \noindent
  for all $x,y$ in $G$.
  If moreover $A$ is a *-algebra we will say that $\letterfors$ is a
*-partial representation if 
  \zitem $(\s x)^* = \si x$ for all $x$ in $G$.  

Observe that under \lcite{\DefPreps.iv} one has that
\lcite{\DefPreps.ii} and \lcite{\DefPreps.iii} become equivalent.
   Given a partial representation $\letterfors$ of $G$ on an algebra
$A$ one has the following useful commutation relation
  $$
  \s x \e y = \e{xy} \s x,
  \eqno{(\seqnumbering)}
  \label\BasicCommutingRelation
  $$
  where $\e y:=\s y\si y$ and $\e{xy}$ is similarly defined
  (see \scite{\Inverse}{2.4} for a proof).

  \definition 
  \label\DefineTPA
  \cite{\TPA,\ExelMisha}.\quad A \stress{twisted partial action} of a group
$G$ on an algebra $A$ is a triple
  $$
  \Theta = \Big(\{D_t\}_{t\in G},\ \{\theta_t\}_{t\in G},\
\{\w{r}{s}\}_{(r,s)\in G\times G}\Big),
  $$
  where, for each $t$ in $G$, $D_t$ is a closed two sided ideal in $A$,
$\theta_t$ is an isomorphism from $D_{t^{-1}}$ onto $D_t$, and for each
$(r,s)$ in $G\times G$, $\w{r}{s}$ is an invertible multiplier of $D_r\cap
D_{rs}$, satisfying the following postulates, for all $r$, $s$ and $t$
in $G$:
  \izitem
  \zitem $D_1 = A$ and $\theta_1$ is the identity automorphism of $A$,
  \zitem $\theta_r(D_{r^{-1}}\cap D_s) = D_r\cap D_{rs}$,
  \zitem $\theta_r(\theta_s(a)) = \w{r}{s}\theta_{rs}(a)\w{r}{s}\inv 
  \for{a\in D_{s^{-1}} \cap D_{s^{-1}r^{-1}}}$,
  \zitem $\w{e}{t}=\w{t}{e}=1$,
  \zitem $\theta_r(a\w{s}{t})\w{r}{st}=\theta_r(a)\w{r}{s}\w{rs}{t} \for{a \in
D_{r^{-1}}\cap D_s\cap D_{st}}$.
  \medskip \noindent
  If moreover $A$ is a *-algebra we will say that the above is a
\stress{*-twisted partial action} if for all $t,r,s\in G$
  \zitem $(D_t)^*=D_t$
  \zitem $\theta_t(a^*) = (\theta_t(a))^*$, for all $a$ in
$D_{t\inv}$,
  \zitem $(\w rs)\inv=(\w rs)^*$.

  \definition 
  \cite{\TPA,\ExelMisha}.\quad 
  Given a twisted partial action, as above, the crossed product algebra,
denoted $A\crossproduct_{\Theta} G$, is defined to be the direct sum  
  $$
  A\crossproduct_{\Theta} G = \bigoplus_{g\in G} D_g,
  $$
  with multiplication
  $$
  (a_g\d g) (a_h \d h) = \theta_g\big(\theta_g\inv(a_g)a_h\big)\w  gh\d{gh},
  $$
  for all $a_g\in D_g$ and $a_h\in D_h$, where we denote by $a_g\d g$
the element $a_g$ viewed in the factor $D_g$ of the above direct sum.

See \cite{\TPA} and \cite{\ExelMisha} for more details, including a
proof of associativity of the above algebra under suitable hypotheses.


\section{Generalities about Hecke pairs}
  Throughout this section $G$ will be a group and $H$ a subgroup.
We will denote by $\lq G H$ (respectively $\rq G H$)
the quotient of $G$ by the equivalence relation according to which
$g_1\sim g_2$ if and only if $g_1\inv g_2\in H$ (respectively $g_1 g_2\inv\in
H$).  
  Thus the equivalence classes relative to $\lq GH$ are the so called
\stress{left cosets} $gH$, for $g\in G$.  Speaking of $\rq GH$ one
similarly has the \stress{right cosets} $Hg$.

We will also consider the equivalence relation according to which the
elements $g_1$ and $g_2$ of $G$ are equivalent when there exist $h,k\in H$
such that $g_1=hg_2k$.  The corresponding \stress{double cosets} therefore
have the form $HgH$, for $g\in G$, and the coset space will be denoted $\dq GHH$.

When $H$ is normal in $G$ then all notions coincide but,  having developed
a bias towards right coset spaces, we will insist in using the notation
$\rq GH$  while most people would prefer to use
$\lq GH$.  Moreover, we will adopt the standard fraction notation for
\stress{right} coset spaces, especially in displayed formulas:

  \definition
  \label \FractionNotation
  If $A$ is a subgroup of a group
$B$ we will let
  $$
  \Rq BA :=  \rq BA.
  $$
  A subset $S$ of $G$ will be called a \stress{family of
representatives} for a coset space (such as the ones above) if there
is exactly one member of $S$ in each equivalence class.

  \state  Proposition
  \relabel \DoubleInRightProp
  Let $H$ be a subgroup of a group $G$.  For every $x\in G$ let
  $$
  \H x = x\inv Hx.
  $$
  Given $x\in G$, let
$S$ be a family of representatives for the coset space
$\rq{H}{(\HH x)}$.  Then 
  $$
  HxH = \dot{\bigcup_{h\in S}} \ Hxh,
  \eqmark{\DoubleInRight}
  $$
  where the symbol \ ``\kern 1.8pt $\dot{\bigcup}$'' \ stands for
disjoint union.  Conversely, if $S$ is any subset of $H$ such that
\lcite{\DoubleInRight} holds then it is a family of representatives for
$\rq{H}{(\HH x)}$.

  \proof The inclusion ``$\supseteq$'' in \lcite{\DoubleInRight} is
obvious so let's prove ``$\subseteq$''.  Given $y\in HxH$ write $y =
k_1xk_2$, with $k_1,k_2\in H$.  By assumption there exists $h\in S$
such that $k_2h\inv\in \HH x$, so that $k_2h\inv = x\inv k x$, for
some $k\in H$.  Therefore
  $$
  y = k_1xk_2 = k_1x (x\inv k x h) = k_1 k x h \in Hxh.
  $$
  In order to prove disjointness suppose that $Hxh = Hxk$, for  $h,k$ in $S$.
Then $xh = \ell xk$ for some $\ell\in H$ whence
  $$
  hk\inv = x\inv\ell x \in \HH x,
  $$
  which implies that $h=k$.
  We leave the converse statement for the reader.
  \proofend

  \definition Let $H$ be a subgroup of a group $G$.  We will say that
$(G,H)$ is a \stress{Hecke pair} if for every $x$ in $G$ one (and
hence all) of the following equivalent conditions hold:
  \izitem
  \zitem $HxH$ is a finite union of right cosets,
  \zitem $\rq H{(\HH x)}$ is finite.

One could as well add two other equivalent conditions to the above,
namely that (iii) $HxH$ is a finite union of left cosets, and (iv) 
$\lq H{(\HH x)}$ is finite; but these will not be used here.

  \definition
  \label\DefineRx
  Let $(G,H)$ be a Hecke pair.  
  \izitem
  \zitem We will denote by $R:G\to{\bf N}$, the function defined by
  $$
  R(x) = \left| \rq H{(\HH x)}\right|
  \for x\in G.
  $$
  \zitem We will denote by $\Delta:G\to{\bf Q}$, the function defined by
  $$
  \Delta(x) = {R(x)\over R(x\inv)}
  \for x\in G.
  $$

By \lcite{\DoubleInRightProp} we have that $R(x)$ is also the number
of right cosets in $HxH$.  It should also be noticed that $R(x\inv)$
is the number of left cosets in $HxH$.  Recall from
\scite{Krieg}{I.3.6} that $\Delta$ is a homomorphism into the additive
group of rational numbers.

From now on we fix a Hecke pair $(G,H)$.

  \definition
  \label \Duality
  Denote by $\F{\rq GH}$ any $F$-vector space having a
basis with as many elements as $\rq GH$.  Fix such a basis and denote
it by
  $$
  \B = \big\{\delta_u : u\in \rq GH\big\}.
  $$
  For each right coset $Hg$ we will denote by $\dpH g$ the linear
functional on $\F{\rq GH}$ given by
  $$
  \<\dH t,\dpH g\> = \left\{
  \matrix{1,& \hbox{if } Ht=Hg,\cr
          0,&\hbox{otherwise,}\hfill } \right.
  $$
  for every $t\in G$,
  where we denote the duality between $\F{\rq GH}$ and its dual space
by $\<\cdot,\cdot\>$, as usual.
  
  \state Proposition
  \label\DefineSx
  Given any $x\in G$ there exists a unique linear operator $\s x$ on $\F{\rq GH}$
such that
  $$
  \s x \calcat {\d {Ht}} = {1 \over |S_x|}\sum_{h\in S_x} \dH{xht}
  \for t\in G,
  \eqmark\DefineHeckeOperator
  $$
  where $S_x$ is any (necessarily finite) family of representatives
for $\rq{H}{(\HH x)}$.
  
  \proof Observing that 
  $$
  HxHt = \dot{\bigcup_{h\in S_x}} \ Hxht,
  $$
  by  \lcite{\DoubleInRight},
  we see that the expression given for $\s x \calcat {\d {Ht}}$ in the
statement is
just the average of the basis elements corresponding to the
right cosets making up $HxHt$.  It is therefore  immediate that $\s x$ is well
defined and does not depend on the choice of $S_x$.
  \proofend

It is clear that $\s h$ is the identity operator for each $h$ in $H$.
In fact this is a special case of the following more general fact:

  \state Proposition
  \label\BiinvarianceOfSx
  For every $x,y\in G$  one has  that
  $
  \s {x} = \s{y}
  $
  if and only if $HxH=HyH$.

  \proof Suppose that $HxH=HyH$, so $y=k_1xk_2$ for some $k_1,k_2\in H$.
Letting  $S_x$ be a family of representatives for $\rq{H}{(\HH x)}$
observe that
  $$
  HyH = Hk_1xk_2H = HxH = 
  \bigcup_{h\in S_x} Hxh =
  \bigcup_{h\in S_x} Hk_1\inv yk_2\inv h =
  \bigcup_{h\in S_x} Hyk_2\inv h,
  $$
  so $k_2\inv S_x$ is a family of representatives for $\rq{H}{(\HH
y)}$.  Therefore  for every $t\in G$,
  $$
  \s{y}\calcat {\d {Ht}} =
  {1\over |S_x|} \sum_{h\in S_x} \dH{yk_2\inv ht} =
  {1\over |S_x|} \sum_{h\in S_x} \dH{k_1xht} =
  {1\over |S_x|} \sum_{h\in S_x} \dH{xht} =
  \s{x}\calcat {\d {Ht}}.
  $$
  Conversely suppose that $\s x = \s y$.  Then, since 
  $
  \<\s x\calcat{\dH{}},\dpH x\> \neq 0,
  $
  one necessarily also has 
  $
  \<\s y\calcat{\dH{}},\dpH x\> \neq 0,
  $
  hence there exists some $k\in S_y$ (a family of representatives for
$\rq{H}{(\HH y)}$), such that $Hyk=Hx$, so that $HyH=HxH$.
  \proofend


  \section{The Hecke algebra}
  Let $(G,H)$ be a Hecke pair, fixed throughout this section.

  \definition The Hecke algebra of the pair $(G,H)$, denoted by
$\HGH$, or simply by $\Hecke$ if the pair $(G,H)$ is
understood, is defined to be the sub-algebra of linear operators on
$\F{\rq GH}$ generated by the set $\{\s x: x\in G\}$.

  For every $g\in G$ denote by $\rho(g)$ the ``right multiplication'' operator on $\F{\rq GH}$
given by
  $$
  \rho(g)\calcat {\d {Ht}} = \d {Htg}
  \for t\in G.
  $$
  It is apparent that $\s x$ commutes with $\rho(g)$ for every
$x$ and $g$ in $G$.  It therefore follows  that each $a\in \Hecke$ 
commutes with every $\rho(g)$.

  \state Proposition
  \label\VectorDetermined
  Let $a,b\in\Hecke$.  If for some $t \in G$ one has that
$a(\dH{t})=b(\dH{t})$, then $a=b$.

  \proof
  For every $s\in G$ one has that
  $$
  a(\dH s) =   
  a\big(\rho(t\inv s)\calcat{\dH t}\big)=   
  \rho(t\inv s)\big(a({\dH t})\big)=
  \rho(t\inv s)\big(b({\dH t})\big)=
  b(\dH s),
  $$
  and hence $a=b$.
  \proofend

  In our next definition we will again make use of the linear functionals
$\dpH g$ introduced in  \lcite{\Duality}.

  \definition 
  \label\DefineFa
  For each $a\in\Hecke$, let $f_a$ be the $F$-valued function on $G$
defined by
  $$
  f_a(t) = \<a(\dH{}),\dpH t\>
  \for t\in G,
  $$
  so that 
  $$
  a(\dH{}) = \sum_{Ht\in\rq GH} f_a(t)\dH t.
  $$

  By \lcite{\VectorDetermined} we see that $a$ is completely
determined by $f_a$.  Observe also that by its very definition, $f_a$
is constant on right cosets.  

  \state Proposition
  \label\FiniteDouble
  For every $a$ in $\Hecke$ one has that $f_a$ is constant on each
double coset.  Moreover $f_a$ is supported in the union of finitely
many such double cosets.

  \proof Given $x\in H$, let $\rho'(x)$ be the dual operator of
$\rho(x)$.  It is immediate to verify that $\rho'(x)\calcat{\dpH t} =
\dpH {tx\inv}$.  For $h,k\in H$ we therefore have
  $$
  f_a(kgh) = 
  \<a(\dH{}),\dpH {kgh}\> =
  \<a(\dH{}),\dpH {gh}\> =
  \<a(\dH{}),\rho'(h\inv)\calcat{\dpH g}\> \$=
  \<\rho(h\inv)\big(a(\dH{})\big),\dpH g\> =
  \<a(\dH{h\inv}),\dpH g\> =
  \<a(\dH{}),\dpH g\> =
  f_a(g).
  $$ 

  As for the last part observe that, since $a(\dH{})$ is a vector in
$\F{\rq GH}$, it is a finite linear combination of the $\dH g$ and
hence $f_a$ is in fact supported in the union of finitely many
\stress{right cosets}, which must obviously involve an even smaller
number of  double cosets.
  \proofend

  We can make use of $f_a$ to describe the matrix of each  operator
$a\in\Hecke$:

  \state Proposition
  \label\TheMatrix
  Let $a\in\Hecke$.  Then, for each $s,t\in G$ one has that 
  $$
  \<a(\dH s),\dpH t\> =
  f_a(ts\inv).
  $$

  \proof
  We have
  $$
  \<a(\dH s),\dpH t\> =
  \<a\big(\rho(s)\calcat{\dH{}}\big),\dpH t\> =
  \<a(\dH{}),\rho'(s)\calcat{\dpH t}\> =
  \<a(\dH{}),\dpH {ts\inv}\> =
  f_a(ts\inv).
  \proofend
  $$

  For the generating operators $\s x$ we have:

  \state Proposition
  \label\CharactOfDouble
  If $x\in G$ then   the function $f_{\s x}$ coincides with the characteristic function
of $HxH$ divided by $|\rq H{(\HH x)}|$.

  \proof Let $S_x$ be a family of representatives for $\rq{H}{(\HH x)}$.
For every $t\in G$ we have
  $$
  f_{\s x}(t) = 
  \<\s x\calcat{\dH{}},\dpH t\> =
  \left\<{1 \over |S_x|}\sum_{h\in S_x} \dH{xh},\dpH t\right\> =
  {1 \over |S_x|} \big[\exists\ h\in S_x,\ \dH{xh}= \dH{t}\big] = \ldots
  $$
  where  the brackets correspond to the boolean value of the logical
statement inside.  Still making use of brackets, the above equals 
  $$
  \ldots = 
  {1 \over |S_x|} \big[Ht\subseteq HxH\big] =
  {1 \over |S_x|} \big[t\in HxH\big].
  $$
  Since $|S_x|$ coincides with $|\rq H{(\HH x)}|$, the proof
is complete.
  \proofend

  We therefore have a description of $H$, at least as far as its
linear structure is concerned:

  \state Proposition
  The correspondence $a\mapsto f_a$ establishes a bijective linear 
correspondence between $\Hecke$ and the space of functions on $G$
which are constant on double cosets and whose support consist of a
finite union of such cosets.

  \proof
  By \lcite{\FiniteDouble} we have that $f_a$ does belong to the
indicated set, while  \lcite{\VectorDetermined} shows that the
correspondence is one-to-one.  That our map is surjective follows from 
\lcite{\CharactOfDouble}.
  \proofend

As an easy consequence we have:

\state Corollary
  \label\LinearBasis
  Let $S$ be a family of representatives for $\dq{G}{H}{H}$. Then the set
$\{\s x: x\in S\}$ is a linear basis for the 
Hecke algebra  $\HGH$.

In order to describe the multiplicative structure of $\Hecke$ in terms of
doubly invariant functions we need the following:

  \state Proposition
  \label\Convolution
  If $a,b\in\Hecke$ then 
  $$
  f_{ab}(t) =
  \sum_{Hs\in\rq GH} f_a(ts\inv) f_b(s)  
  \for t\in G.
  $$

  \proof
  We have 
  $$
  f_{ab}(t) =
  \<ab(\dH{}),\dpH t\> =
  \Big\<a\Big(\sum_{Hs\in\rq GH} f_b(s)\dH s\Big),\dpH t\Big\> \$=
  \sum_{Hs\in\rq GH} f_b(s)\<a(\dH s),\dpH t\> \={(\TheMatrix)}
  \sum_{Hs\in\rq GH}   f_a(ts\inv) f_b(s),
  $$
  concluding the proof.
  \proofend

We thus reconcile our point of view with the classical definition of
Hecke algebras (see e.g.~\scite{Krieg}{I.4}):

  \state Corollary
  \label\HeckeAsDoublyInvFunctions
  $\Hecke$ is isomorphic to the algebra of doubly invariant functions
on $G$ which are supported in the union of finitely many
double cosets, equipped with the convolution product defined, for
every $f$ and $g$ in said algebra, by
  $$
  (f * g)(t) =
  \sum_{Hs\in\rq GH} f(ts\inv) g(s)  
  \for t\in G.
  $$


  \section{$*$-algebra structure} 
  In this section we will turn  the Hecke algebra  into a *-algebra.
The reader might be familiar with the involution given by 
  $$
  f^\#(t) = \conj {f(t\inv)},
  \eqno{(\seqnumbering)}
  \label\UsualStar
  $$
  used by many authors (see e.g. \cite{\BostConnes}).
  However the involution used in \cite{\QuiggKaliLand} is better
suited  for our purposes, given our emphasis on the operators $\s x$,
as we shall see shortly.

Below we will use the rational homomorphism $\Delta(x) =
R(x)/R(x\inv)$ defined in \lcite{\DefineRx.ii}.

  \definition 
  \label\NewInnerProduct
  We will denote by $\ip[\cdot,\cdot]$ the unique sesqui-linear form
on $\F{\rq GH}$ such that for every $t,s\in G$,
  $$
  \ip[\dH t,\dH s] = 
  \left\{
  \matrix{
    \vrule depth  20pt width 0pt
    \Delta(s)&\hbox{, if } Ht = Hs, \cr
    0&\hbox{, if } Ht\neq Hs. \hfill
     } \right.
  $$

  Observe that when $Ht=Hs$, then we have that $HtH=HsH$ so
  $R(t)/R(t\inv)=  R(s)/R(s\inv)$
  and we see that  our  form is hermitian.
  It is elementary to verify that it is non-degenerate
as well.

  \state Proposition
  \label\Adjoint
  For every $x\in G$ one has
  $$
  \ip[\s x (\xi),\eta] =
  \ip[\xi, \s {x\inv}(\eta)]
  \for \xi,\eta\in \F{\rq GH}.
  $$

  \proof It is obviously enough to consider $\xi=\dH t$, and $\eta =
\dH s$, where $t,s\in G$.

  Choose families of representatives $S_x$ and $S_{x\inv}$ for the
coset spaces $\rq{H}{(\HH x)}$ and $\rq{H}{(\HH {x\inv})}$, respectively,
so that $|S_x| = R(x)$ and $|S_{x\inv}|=R(x\inv)$.
  We have
  $$
  \ip[\s x\calcat{\dH t},\dH s] =
  {1 \over R(x)}\sum_{h\in S_x} \ip[\dH{xht},\dH s] =
  {R(s) \over R(x)R(s\inv)}\ [st\inv\in HxH],
  $$
  where the brackets denote boolean value, as before.  On the other
hand
  $$
  \ip[\dH t,\s{x\inv}\calcat{\dH s}] =
  {1 \over R(x\inv)}\sum_{k\in S_{x\inv}}\ \ip[\dH t,\dH {x\inv ks}] =
  {R(t) \over R(x\inv) R(t\inv)}\ [ts\inv\in Hx\inv H] \$=
  {R(t) \over R(x\inv) R(t\inv)}\ [st\inv\in Hx H].
  $$

  In order to complete the proof it is then enough to prove that 
  $$
  {R(s) \over R(x)R(s\inv)} =
  {R(t) \over R(x\inv) R(t\inv)}
  \eqmark\AdjointBoilsDown
  $$
  whenever $st\inv\in HxH$.

  Write $st\inv=hxk$, with $h,k\in H$, so that
$x=h\inv st\inv k\inv$ and, given that $R$ is clearly a doubly
invariant function we have that $R(x) = R(h\inv st\inv k\inv) =
R(st\inv)$ and hence  \lcite{\AdjointBoilsDown} boils down to
  $$
  {R(s) \over R(st\inv)R(s\inv)} =
  {R(t) \over R(ts\inv) R(t\inv)}
  $$
  which follows immediately from \scite{Krieg}{I.3.6}.
  \proofend

  \state Corollary 
  \label\StarAlgebraStructure
  For every $a$ in $\HGH$ there exists a unique $a^*$ in $\HGH$ such
that
  $$
  \ip[a(\xi),\eta] =
  \ip[\xi, a^*(\eta)]
  \for \xi,\eta\in \F{\rq GH}.
  $$
  In addition $\HGH$ becomes a *-algebra under the operation $a\mapsto
a^*$ and,  for every $x\in G$, we have  $\s x^*=\si x$.

With this we may improve the description of
$\HGH$ in terms of doubly invariant functions given in
\lcite{\HeckeAsDoublyInvFunctions}:

  \state Proposition
  \label\NewStar
  $\HGH$ is *-isomorphic to the algebra of doubly invariant
functions described in \lcite{\HeckeAsDoublyInvFunctions} once the
latter  is
made a *-algebra by the involution  given by 
  $$
  f^*(t) = \Delta(t\inv)\ \conj{f(t\inv)}
  \for t\in G,
  $$
  for every $f$ in said function algebra.

  \proof  Observe that for all $t$ and $s$ in $G$ we have that
  $$
  \ip[\dH t,\dH s] = \Delta(s)\ \<\dH t,\dpH s\>,
  $$
  where the duality in the right-hand-side is given by \lcite{\Duality}.
  So obviously 
  $$
  \ip[\xi,\dH s] = \Delta(s)\ \<\xi,\dpH s\>
  \for \xi\in  \F{\rq GH}.
  $$
  Given $a$ in $\HGH$ we have by definition \lcite{\DefineFa} that
  $$
  f_{a^*}(t) =
  \<a^*(\dH{}),\dpH t\> = 
  \Delta(t\inv)\ \ip[a^*(\dH{}),\dH t] =
  \Delta(t\inv)\ \conj{\ip[a(\dH t),\dH{}]} \$=
  \Delta(t\inv)\ \conj{\<a(\dH t),\dpH{}\>} \ \={(\TheMatrix)}\
  \Delta(t\inv)\ \conj{f_a(t\inv)} =
  (f_a)^*(t).
  \proofend
  $$

Our last result of this section shows that, under certain hypotheses
about the field $F$, the two involutions are essentially the same:

  \state Proposition
  \relabel\CompareStars
  Suppose there  exists a group
homomorphism $\lambda$ from $G$ to the multiplicative group
of $F$ such that $\lambda(x)^2=\Delta(x)$, for all $x$ in
$G$.  Then the *-algebras 
  $\Big(\HGH, * \Big)$ and $\Big(\HGH, \#\Big)$
  are isomorphic.

  \proof
  It is elementary to check that  the map
  $$
  \Lambda : \Big(\HGH, * \Big) \to \Big(\HGH, \#\Big)
  $$
  given by 
  $$
  \Lambda(f)\Calcat x = \lambda(x) f(x)
  \for x\in G,
  $$
  for all $f$ in $\HGH$, is an isomorphism of *-algebras.
  \proofend


  \section{Commuting subgroups}
  \label\CommutingGpsSection
  We will now develop a few basic facts about commuting subgroups in
preparation for our study of protonormal subgroups.

\definition 
  \izitem 
  \zitem If $A$ and $B$ are subsets of a group $G$ we will denote
by $AB$ the set
  $$
  AB =\{ab: a\in A,\ b\in B\}.
  $$
  \zitem If $A$ and $B$ are subgroups of $G$ we will say that $A$ and
$B$ \stress{commute} if $AB = BA$.

\bigskip The following lists useful alternative  characterizations of the concept
above:

\state Proposition 
  \label\AltCommutSGps
  Given subgroups $A$ and $B$ of a group $G$ the following are
equivalent
  \izitem
  \zitem $A$ and $B$ commute,
  \zitem $BA\subseteq AB$,
  \zitem $AB$ is closed under multiplication,
  \zitem $AB$ is a subgroup of $G$.
  
  \proof

  (i)$\Rightarrow$(ii): obvious.

  (ii)$\Rightarrow$(iii): we have
  $$
  AB AB = A(BA)B \subseteq A(AB)B = AABB = AB.
  $$

  (iii)$\Rightarrow$(iv): we have
  $$
  (AB)\inv = B\inv A\inv = BA \subseteq ABAB \subseteq AB,
  $$
  so $AB$ is closed under taking inverses and hence is a subgroup.

  (iv)$\Rightarrow$(i).  We have
  $$
  BA\subseteq ABAB = AB.
  $$
  Taking inverses we get $AB\subseteq BA$, so $AB=BA$.
  \proofend


We now list two elementary results for future reference, in which the
fraction notation introduced in \lcite{\FractionNotation} is used.

  \state Lemma
  \label\TwoQuotients
If the subgroups $A$ and $B$ commute there is a natural bijection
  $$
  \Rq{B}{A\cap B} \to 
  \Rq{AB}{A}
  $$
  which sends the right coset $(A\cap B)b$ to the right coset $Ab$,
for every $b\in B$.

\proof Left to the reader. \proofend

  \state Lemma
  \label\ThreeGroups
  Let $A$, $B$, and $C$ be groups with $A\subseteq B\subseteq C$ and
let $\{b_i: i\in I\}$ and $\{c_j: j\in J\}$ be families of
representatives for the coset spaces $\rq BA$ and $\rq CB$,
respectively.  Then $\{b_ic_j: (i,j)\in I\times J\}$ is a family of
representatives for $\rq CA$.  In particular, if $\rq CA$ is finite,
then $\rq BA$ and $\rq CB$ are both finite and
  $$\left|\Rq CA\right| = \left|\Rq BA\right|\left|\Rq CB\right|.$$

\proof Left to the reader. \proofend  

Let us fix, for the time being, a group $G$ and a subgroup $H$ and
let $\F{\rq GH}$ be as defined in \lcite{\Duality}.

\definition If $S$ is any finite subset of $\rq GH$ we will denote by
$\mu(S)$ the average of the elements of $S$ computed in $\F{\rq GH}$.  Precisely speaking,
  $$
  \mu(S) = {1\over |S|} \sum_{s\in S} \delta_s.
  $$

\definition Given a subgroup $K$ of $G$ which commutes with $H$ and
such that $\rq K{(H\cap K)}$ is finite, observe that $\rq{HK}{H}$ is a
finite subset of $\rq GH$ by \lcite{\TwoQuotients}.  We therefore  denote by
$\q{K}$ the element of $\F{\rq GH}$ defined by
  $$
  \q{K} = \mu(\rq{HK}{H}).
  $$

  If $S\subseteq K$ is a family of representatives for the coset space
$\rq{K}{(H\cap K)}$, then by \lcite{\TwoQuotients} we have that the
elements of the form $Hk$, with $k\in S$, are precisely all of the
(pairwise 
distinct) elements of $\rq{HK}{H}$ and hence 
  $$
  \q{K} = {1\over|S|}\sum_{k\in S}\dH k.
  \eqmark\Alternateq
  $$

In addition to $H$ we will now fix a subgroup $K$ of $G$ as above,
that is, such that $K$ commutes with $H$ and $\rq K{(H\cap K)}$ is
finite.

As before let us denote by $\rho$ the right-regular
(anti-)representation of $G$ on $\F{\rq GH}$.

\state Proposition
  \label\TocomPressa
  For all $g$ in $HK$ one has that 
  $$
  \rho_g(\q{K}) = \q{K}.
  $$

  \proof If $g\in HK$ then the operator $\rho_g$ clearly leaves $\rq
{HK}H$ invariant and hence it must consist of a permutation of the
elements in the latter set, therefore leaving $\q{K}$ unchanged.
  \proofend

If $x,y\in G$ are such that $Hx=Hy$, then $x=hy$ for some $h\in H$ and
hence
  $$
  \rho_{x}(\q{K}) =
  \rho_{hy}(\q{K}) =
  \rho_y\big(\rho_{h}(\q{K})\big) \={(\TocomPressa)}
  \rho_{y}(\q{K}),
  $$
  so the expression $\rho_{x}(\q{K})$ depends only on the right coset
where $x$ lies.  This proves the following:

\state Proposition  
  The correspondence
  $$
  x\in G \longmapsto \rho_{x}(\q{K})\in \F{\rq GH}
  $$
  drops to the quotient providing a well defined map from $\rq GH$ to
$\F{\rq GH}$ which, when linearized, gives an operator $\Q{K}$ on
$\F{\rq GH}$ satisfying
  $$
  \Q{K}(\dH x) = \rho_x(\q{K}).
  \eqmark{\CharacterizeQ}
  $$

If $S\subseteq K$ is a family of representatives for the coset space
$\rq{K}{(H\cap K)}$ as in \lcite{\Alternateq}, notice that for all
$x\in G$,
  $$
  \Q{K}(\dH x) =
  \rho_x(\q{K}) =
  {1\over|S|}\sum_{k\in S}\rho_x(\dH k) =
  {1\over|S|}\sum_{k\in S}\dH {kx}.
  \eqmark\qWithReps
  $$

  This should be compared to the identity
  $
  \big(\sum_{k\in S}\dH k\big)\dH x =
  \sum_{k\in S}\dH {kx},
  $
  which would only make sense if $H$ were a normal subgroup of $G$ and we
were  using of the group-algebra structure of $\F{\rq GH}$.

Denote by
  $
  \pi:\rq G{H} \to \rq G{HK},
  $
  the quotient map and let 
  $$
  \pit:\F{\rq G{H}} \to \F{\rq G{HK}}
  $$
  be its linearization.

  \state Proposition
  \label\IsoWhenRestricted
  The restriction of \ $\pit$ to the range of\/ $\Q{K}$ is a linear
isomorphism onto $\F{\rq G{HK}}$. In addition
  $$
  \pit(\Q{K}(\dH x)) = \delta_{HKx}
  \for x\in G,
  $$
  so that $\pit \circ \Q K = \pit$.

  \proof
  Denote by $\hat\pi$ the restriction of\/ $\pit$ to the range of\/
$\Q{K}$.
  Let $S$ be as in \lcite{\qWithReps} so that for all $x$ in
$G$ one has
  $$
  \hat\pi(\Q{K}(\dH x)) = 
  {1\over|S|}\sum_{k\in S}\pit(\dH {kx}) =
  {1\over|S|}\sum_{k\in S}\delta_{HKkx} =
  \delta_{HKx},
  $$
  where the last step holds because $S\subseteq K$.   This proves
the identity in the statement and also that $\hat\pi$ is surjective.  In
order to prove injectivity consider the map
  $$
  \phi : x\in G \longmapsto \Q{K}(\dH x)\in  \F{\rq GH},
  $$
  and observe that if $x,y\in G$ are such that $x = gy$, with $g\in HK$, then
  $$
  \phi(x)=
  \Q{K}(\dH x) =
  \rho_{x}(\q{K}) =
  \rho_{gy}(\q{K}) =
  \rho_y\big(\rho_g(\q{K})\big) \={(\TocomPressa)}
  \rho_y(\q{K}) =
  \Q{K}(\dH y) = 
  \phi(y).
  $$
  Therefore $\phi$ drops to the quotient $\rq G{HK}$ and the
corresponding linearization is a map
  $$
  \tilde\phi:\F{\rq G{HK}} \to \F{\rq GH}
  $$
  satisfying
  $$
  \tilde\phi(\delta_{HKx}) = \Q{K}(\dH x)
  \for x\in G.
  $$
  Therefore we have  for all $x\in G$, that
  $$
  \tilde\phi\Big(\hat\pi\big(\Q{K}(\dH x)\big)\Big) =
  \tilde\phi(\delta_{HKx}) =
  \Q{K}(\dH x)
  $$
  showing that $\tilde\phi \circ \hat\pi$ is the identity map on the
range of $\Q{K}$.  Thus $\hat\pi$ is injective.
  \proofend

In the last result of this section we shall again refer to the
quotient map
  $
  \pi:\rq G{H} \to \rq G{HK},
  $
  as well as to its linearized version $\pit$.

  \state Proposition 
  \label\ThreeCommutingGroups
  If $L$ is yet another subgroup of\/ $G$ which commutes with both $H$
and $K$, and such that $\rq L{(H\cap L)}$ is finite, then 
  $$
  \pit(\q{L}) = \mu(\rq{HKL}{HK}). 
  \eqmark\ThreeCommutingGroupsIdentity
  $$
  In particular $\pit(\q{L})$ only
depends on the image of $\rq{HL}H$ under $\pi$.

\proof Consider the chain of subgroups
  $$
  H\cap L \ \subseteq\  HK\cap L \ \subseteq\  L
  $$
  and let $\{b_i: i\in I\}$ and $\{c_j: j\in J\}$ be families of
representatives for the coset spaces $\rq {(HK\cap L)}{(H\cap L)}$ and
$\rq {L}{(HK\cap L)}$, respectively.   By \lcite{\ThreeGroups} we then
have that 
$\{b_ic_j: (i,j)\in I\times J\}$ is a family of representatives for
$\rq L{(H\cap L)}$, which is a finite set by hypothesis hence implying
that both $I$ and $J$ must be finite sets as well.  By
\lcite{\Alternateq} we have that
  $$
  \pit(\q{L}) =
  \pit\left({1\over|I||J|}\sum_{i\in I}\sum_{j\in J}\dH{b_ic_j}\right)
=
  {1\over|I||J|}\sum_{i\in I}\sum_{j\in J}\delta_{HKb_ic_j} \$=
  {1\over|J|}\sum_{j\in J}\delta_{HKc_j} =
  \mu(\rq{HKL}{HK}),
  $$
  where the last step follows from the natural equivalence between
$\rq {L}{(HK\cap L)}$ and $\rq{HKL}{HK}$ given by
\lcite{\TwoQuotients}.
  \proofend


 \section{Protonormal subgroups}

  \definition
  \label \ProtoNormal
  Let $H$ be a subgroup of a group $G$.  We will say that $H$ is a
  \stress{protonormal}\fn{From ``{\tt dictionary.reference.com}'':
proto- ({\it pref.}) 4. Having the least amount of a specified element
or radical.  }
  subgroup if $\H x$ and $H$ commute for every $x\in G$ (recall from
\lcite{\DoubleInRightProp} that $\H x$ means $x\inv H x$).

Observe that every normal subgroup $H$ is protonormal since $\H x = H$
for all $x$ in $G$.  More generaly, suppose that there exists a
subgroup $N$ of $G$ containing $H$ such that $H\normal N\normal G$
(the symbol ``$\,\normal$'' standing for ``is normal in''), in which
case it is sometimes customary to say that $H$ is 2-subnormal, which
we shall shorten to subnormal.  Then for every $x$ in $G$ and $h$ in
$H$ we have that $x\inv hx\in N$ and hence 
  $$
  x\inv hx H =
  Hx\inv hx.
  $$
  It easily follows  that $\H x$ and $H$ commute.  In other words,
every subnormal subgroup is protonormal.

Given $y\in G$ and assuming that $\H{yx\inv}$ and $H$ commute  we
conclude, upon applying the inner automorphism,
  $$
  {Ad}_{x\inv}: g \in G\mapsto  x\inv g x \in G,
  $$
  that $\H y$ and $\H x$ also commute.  Thus, if $H$ is a protonormal
subgroup then  all of its conjugates commute among themselves.  It is also
evident that the subgroups of the form
  $$
  \H{x_1} \H{x_2}\ldots \H{x_n},
  $$
  where $x_1,x_2,\ldots,x_n\in G$, all commute with each other.
  Since for all $y\in G$ we have that
  $$
  y\inv(\H{x_1} \H{x_2}\ldots \H{x_n})y = 
  \H{x_1y} \H{x_2y}\ldots \H{x_ny},
  $$
  we see that $\H{x_1} \H{x_2}\ldots \H{x_n}$ is also protonormal.

When $(G,H)$ is a Hecke pair such that $H$ is protonormal in $G$ we
have by definition that $\rq H{(\HH x)}$ is finite and therefore so is
$\rq{\H xH}{\H x}$ by \lcite{\TwoQuotients}.  This allows for a
slightly different but usefull description for the  operators $\s
x$ of \lcite{\DefineSx}:

  \state Proposition
  \label\AlternateHeckeOperator
  If\/ $T_x$ is a family of representatives for the coset space $\rq{\H
xH}{\H x}$, then 
  $$
  \s x \calcat {\d {Ht}} = {1 \over |T_x|}\sum_{k\in T_x} \dH{xkt}
  \for t\in G.
  $$

  \proof
  Let $T_x = \{k_1,\ldots,k_n\}$ and write each $k_i$ as $\ell_i h_i$,
with $\ell_i\in\H x$ and $h_i\in H$.  It is then easy to prove  that
$\{h_1,\ldots,h_n\}$ is a family of representatives for
$\rq{H}{(\HH x)}$.  In addition notice that $x\ell_ix\inv\in H$, so
that
  $$
  Hxk_it =
  Hx\ell_ih_it =
  Hx\ell_ix\inv xh_it =
  Hxh_it,
  $$
  from where the result follows.
  \proofend


  \section{The canonical partial representation}  
  Throughout this section we will fix a Hecke pair $(G,H)$ such that
$H$ is a protonormal subgroup of $G$.  Our major goal will be to
show that $\letterfors$ is a partial group representation.

  \state Theorem
  \relabel\MainPrep
  If $(G,H)$ is a Hecke pair with $H$ protonormal in $G$ then the
correspondence
  $$
  x\in G \mapsto \s x \in \HGH
  $$
  is a partial representation.

  \proof 
  Axiom \lcite{\DefPreps.i} is obviously verified so we begin by
proving that for every $x$ and $y$ in $G$ one has that
  $$
  \s{x\inv}\s x\s y = \s{x\inv}\s {xy}.
  $$
  By \lcite{\VectorDetermined} it is enough to show that these
operators coincide on $\dH {y\inv}$.  

For every $u\in\{x\inv, x,y,xy\}$, pick a family of representatives
$S_u$ for the coset space $\rq{H}{(\HH u)}$.
  We therefore have
  $$
  \s{x\inv}\s x\s y \calcat{\dH {y\inv}} =  
  {1\over |S_{x\inv}|\ |S_x|\ |S_y|}
  \sum_{h\in S_{x\inv}} \sum_{k\in S_x} \sum_{\ell\in S_y}
  \dH{x\inv hxky\ell y\inv} = \ldots
  $$
  Recalling that $\rho$ denotes the right regular representation of $G$ on
$\F{\rq GH}$ we may write  the above as
  $$
  \ldots =
  {1\over |S_x|\ |S_y|}
  \sum_{\ell\in S_y} \sum_{k\in S_x} \rho_{ky\ell y\inv}\left(
  {1\over |S_{x\inv}|} \sum_{h\in S_{x\inv}} 
  \dH{x\inv hx}\right) = \ldots
  $$
   Given that $S_{x\inv}$ is a family of representatives for 
$\rq{H}{(\HH {x\inv})}$, it is evident that 
  $\{x\inv hx : h\in S_{x\inv}\}$ is a family of representatives for
$\rq{\H x}{(\HH {x})}$, so that the term within the big pair of
parenthesis above coincides with $\q{\H x}$ by \lcite{\Alternateq}.
Here we are using the results of section \lcite{\CommutingGpsSection}
with  the role of the  groups $H$ and $K$ mentioned there played by $H$ and $\H
x$, respectively.
The above then equals
  $$
  \ldots =
  {1\over |S_x|\ |S_y|}
  \sum_{\ell\in S_y} \sum_{k\in S_x} \rho_{ky\ell y\inv}(\q{\H x}) = 
  {1\over |S_x|\ |S_y|}
  \sum_{\ell\in S_y} \sum_{k\in S_x} \rho_{y\ell y\inv}(\rho_k(\q{\H
x})) \= {(\TocomPressa)} $$ $$ =
  {1\over |S_y|} \sum_{\ell\in S_y} \rho_{y\ell y\inv}(\q{\H x}) = 
  \Q{\H x}\left({1\over |S_y|} \sum_{\ell\in S_y} \dH{y\ell
y\inv}\right) =
  \Q{\H x}\left(\q{\H{y\inv}}\right),
  $$
  where the last identity again follows from \lcite{\Alternateq} since $\{y
\ell y\inv : \ell\in S_{y}\}$ is a family of representatives for
  $\rq{\H {y\inv}}{(\HH {y\inv })}$.

  On the other hand
  $$
  \s{x\inv}\s {xy} \calcat{\dH {y\inv}} =  
  {1\over |S_{x\inv}|\  |S_{xy}|}
  \sum_{h\in S_{x\inv}} \sum_{m\in S_{xy}}
  \dH{x\inv hxymy\inv} \$=
  {1\over |S_{xy}|}  \sum_{m\in S_{xy}}  
  \rho_{ymy\inv} \left(
  {1\over |S_{x\inv}|} \sum_{h\in S_{x\inv}} 
  \dH{x\inv hx}\right) =
  {1\over |S_{xy}|}  \sum_{m\in S_{xy}}  
  \rho_{ymy\inv}(\q{\H x}) \$=
  \Q{\H x}\left( {1\over |S_{xy}|} \sum_{m\in S_{xy}} \dH{ymy\inv}
\right).
  $$

  \def\lostq{q'}
  Denoting by $\lostq$ the element of $\F{\rq GH}$ enclosed by the last big pair of parenthesis
above our task is therefore reduced to proving the identity
  $$
  \Q{\H x}\left(\q{\H{y\inv}}\right) = \Q{\H x}(\lostq).
  $$
  Employing \lcite{\IsoWhenRestricted} we see that the above identity
holds if and only if
  $
  \pit\left(\Q{\H x}\left(\q{\H{y\inv}}\right)\right) = \pit(\Q{\H x}(\lostq)),
  $
  which is to say that
  $$
  \pit\left(\q{\H{y\inv}}\right) = \pit(\lostq),
  \eqmark\MainGoal
  $$
  by the last part of \lcite{\IsoWhenRestricted}.


  Consider the diagram below in which we use the notation described in
\lcite{\FractionNotation}:

  \def\frac#1#2{\ds{#1\over #2}}
  $$
  \def\spc{\ }
  \def\rght{\spc\longrightarrow\spc}
  \def\lrght#1{\spc{\buildrel #1\vrule depth 5pt width 0pt \over\longrightarrow}\spc}
  \matrix{
  \frac{H}{\HH{y}} & \lrght{Ad_{y}} &
  \frac{\H{y\inv}}{\H{y\inv}\cap H} & \lrght{(\TwoQuotients)} &
  \frac{H\H{y\inv}}{H}&\ \hookrightarrow &
  \frac{G}{H} \cr\cr
  &&&&&&\ \big\downarrow\pi \cr\cr
  &&&&&&\frac{G}{H\H x} \cr\cr
  &&&&&&\ \ \big\uparrow \pi_x \cr\cr
  \frac{H}{\HH{xy}} & \lrght{Ad_{y}} &
  \frac{\H{y\inv}}{\H{y\inv}\cap\H{x}} & \lrght{(\TwoQuotients)} &
  \frac{\H x\H{y\inv}}{\H x}&\ \hookrightarrow &
  \frac{G}{\H x} \cr\cr
  }
  $$
  where the arrows ``$\hookrightarrow$'' refer to inclusion and the
vertical arrows are quotient mappings.  

We now intend to apply \lcite{\ThreeCommutingGroups} for the two
situations outlined in the rows in our diagram.  
Precisely, with respect to the top row, the triple $(H,K,L)$ of groups referred
to in \lcite{\ThreeCommutingGroups} will be taken to be
  $(H,\H x,\H{y\inv})$.
  Identity \lcite{\ThreeCommutingGroupsIdentity} is then translated to
  $$
  \pit\left(\q{\H{y\inv}}\right) = \mu(\rq{H\H x\H{y\inv}}{H\H x}).
  $$

Speaking of the bottom row, take the triple $(H,K,L)$ of
\lcite{\ThreeCommutingGroups} to be $(\H x,H,\H{y\inv})$.   In order
to distinguish from the previous application of
\lcite{\ThreeCommutingGroups}, we will use $\qalt{}$ in place of
$\q{}$.

Observe that since $S_{xy}$ is a family of
representatives for $\rq{H}{(\HH{xy})}$, we have that $Ad_{y}(S_{xy})$
is a family of representatives for
$\rq{\H{y\inv}}{(\H{x}\cap\H{y\inv})}$.  Therefore
  $$
  \qalt{\H{y\inv}} =
  {1\over |S_{xy}|} \sum_{m\in S_{xy}} \delta_{\H xymy\inv}.
  $$
  Applying \lcite{\ThreeCommutingGroups} we therefore deduce that
  $$
  \pit_x\left(\qalt{\H{y\inv}}\right) =  \mu(\rq{\H xH\H{y\inv}}{\H xH}).
  $$
  Since the $\mu$'s of our two situations coincide, as they both
correspond to averaging within $\F{\rq G{H\H x}}$, we then conclude
that
  $$
  \pit\left(\q{\H{y\inv}}\right) =  \pit_x\left(\qalt{\H{y\inv}}\right).
  $$
  It follows that
  $$
  \pit\left(\q{\H{y\inv}}\right) =  
  \pit_x\left(\qalt{\H{y\inv}}\right) =
  {1\over |S_{xy}|} \sum_{m\in S_{xy}} \pit_x(\delta_{\H xymy\inv}) \$=
  {1\over |S_{xy}|} \sum_{m\in S_{xy}} \delta_{H\H xymy\inv} =
  {1\over |S_{xy}|} \sum_{m\in S_{xy}} \pit(\delta_{H ymy\inv}) =
  \pit(\lostq),
  $$
  proving \lcite{\MainGoal} and hence showing that $\letterfors$
satisfies \lcite{\DefPreps.ii}.  With respect to
\lcite{\DefPreps.iii} observe that for all $\xi,\eta\in\F{\rq GH}$ one
has
  $$
  \ip[\s x\s y \s{y\inv}(\xi),(\eta)] \={(\Adjoint)}
  \ip[\xi,\s y\s{y\inv}\s{x\inv}(\eta)] =
  \ip[\xi,\s y\s{y\inv x\inv}(\eta)] \={(\Adjoint)}
  \ip[\s {xy} \s{y\inv}(\xi),\eta].
  $$ 
  Given that $\ip[\cdot,\cdot]$ is nondegenerated we conclude that
  $$
  \s x\s y \s{y\inv} =
  \s {xy} \s{y\inv}.
  \proofend
  $$

We would now like to show that it is necessary to assume that $H$ is
protonormal in $G$ in order to conclude that $\letterfors$ is a
partial representation.

\state Theorem
  \relabel\NecessaryCondForPrep
  Let $(G,H)$ be a Hecke pair such that
  $$
  \s x\s{x\inv}\s x = \s x
  \for x\in G,
  $$
  (which obviously holds in case $\letterfors$ is a partial
representation).  Then $H$ is protonormal in $G$.

  \proof  For every $u\in\{x\inv, x\}$, pick a family of representatives
$S_u$ for the coset space $\rq{H}{(\HH u)}$.  Then 
  $$
  \s x\s{x\inv}\s x \calcat{\dH {}} =  
  {1\over |S_x|^2|S_{x\inv}|}
  \sum_{h\in S_x} \sum_{k\in S_{x\inv}} \sum_{\ell\in S_x}
  \dH{x hx\inv kx\ell}.
  $$
  It is easy to see that every right coset contained in $HxHx\inv HxH$
occurs with a nonzero coefficient in the sum above\fn{Observe that
we are using, in a non-trivial manner, that the characteristic of
$\Fd$ is zero.}.  It must
therefore occur as well in the sum describing $\s x\calcat {\dH{}}$, namely
  $$
  \s x \calcat {\dH{}} = {1 \over |S_x|}\sum_{h\in S_x} \dH{xh}.
  $$
  It follows that
  $HxHx\inv HxH \subseteq HxH$.  Multiplying this on the left by
$x\inv$ gives
  $$
  x\inv HxHx\inv HxH \subseteq x\inv HxH,
  $$
  or equivalently that
  $\H xH\H xH \subseteq \H xH$.  Using \lcite{\AltCommutSGps.iii} it follows that
$\H x$ and $H$ commute, so $H$ is protonormal in $G$, as desired.
  \proofend

Observe that the kernel of $\letterfors$, namely
  $$
  {\rm Ker}(\letterfors) = \{x\in G: \s x = 1\}
  $$
  is precisely $H$.  This shows that, while the kernel of a partial
representation is always a subgroup, it needs not be normal.  This
motivates the general question as to which subgroups of a group $G$
coincide with the kernel of a partial representation.  The answer is
very simple, all subgroups do.  Given any subgroup $H\subseteq G$
consider the map $u:G\to F$ given by
  $$
  u(x) = \left\{
  \matrix{1,& \hbox{if } x\in H,\cr
          0,&\hbox{otherwise.}\hfill } \right.
  $$
  It is easy to see that $u$ is a partial representation and clearly 
  ${\rm Ker}(u) = H$.


\section{Generalities about subnormal groups}
  Some of our results can only be proved for subgroups which are a bit
more normal than protonormal.  We shall briefly describe this class in
what follows referring the reader to \cite{\Ward} for more information.

  \definition
  \label \DefSubnormal
  Let $H$ be a subgroup of a group $G$.  We will
say that $H$ is \stress{subnormal} in $G$ if for every $x\in G$ and
$h,k\in H$ one has that
  $$
  xhx\inv k  xh\inv x\inv \in H.
  $$

Writing the above as $(xhx\inv) k (xh x\inv)\inv$, this
says that $H$ is closed under conjugation by elements $g$ in $G$ of
the form $g=xhx\inv$ (which itself is the conjugation of the element
$h\in H$ by the arbitrary element $x\in G$).

  \state Proposition
  \label\AlternateSubnor
  If $H$ is a subgroup of a group $G$ then the following are
equivalent:
  \izitem
  \zitem $H$ is subnormal in $G$.
  \zitem For
every $x\in G$ and $h\in H$ one has that
  $
  Hxhx\inv = xhx\inv H.
  $
  \zitem $H$ is normal in the  intersection of all normal subgroups of\/ $G$
containing $H$.
  \zitem There exists a subgroup $N$ of $G$ such that 
  $H\normal N\normal G$.

  \proof Observe that any normal subgroup of $G$ containing $H$ must
contain the set
  $$
  Y = \{xhx\inv: x\in G,\ h\in H\},
  $$
  and hence also the subgroup  $N$  generated by $Y$.

Since $Y$ is obviously invariant under conjugation by elements of $G$,
one sees that the same applies to $N$, that is, $N\normal G$.
This said it becomes clear that $N$ is the intersection of all normal
subgroups of $G$ containing $H$ mentioned in (ii).  

Assuming (i) notice that $yH y\inv = H$, for every $y$ in $Y$.
Therefore the same holds for every $y$ in $N$.  So $H\normal N$.  This
proves that (i) $\Rightarrow$ (iii).

It is obvious that (iii) $\Rightarrow$ (iv).  In order to
show that (iv) $\Rightarrow$ (i) let $N$ be as in (iv) and let $x\in G$ and
$h,k\in H$.  Observe that the element $n = xhx\inv$ satisfies
  $$
  n=  xhx\inv \in xHx\inv \subseteq xNx\inv = N,
  $$
  so that
  $$
  (xhx\inv) k (xh x\inv)\inv  = nkn\inv \in nHn\inv\subseteq H,
  $$
  because $H\normal N$.

We leave the elementary implication (i) $\Leftrightarrow$ (ii) for the reader.
  \proofend

Recall from the introduction that $H\subseteq G$ is said to be
$n$-subnormal if there exists a normal chain
  $$
  H=H_0\normal H_1\normal H_2\normal\ldots \normal H_n=G,
  $$
  of length $n$.
  Thus, our concept of subnormality is equivalent to 2-subnormality.
Also observe that every subnormal subgroup is protonormal.

  \state Proposition
  \label\HxNormal
  If\/ $H$ is subnormal in $G$ let $N$ be a subgroup of $G$ such that
$H\normal N\normal G$.  Then for every $x\in G$ one has that $\H
x\normal N$.  In particular
  \izitem
  \zitem $\HH x\normal H$, and
  \zitem $H\normal H\H x$.

\proof  $\H x$ is contained in $N$ because
  $$
  \H x = x\inv Hx \subseteq x\inv Nx = N.
  $$
  Moreover $\H x$ is the image of $H$ under the (not necessarily
internal) automorphism $Ad_{x\inv}$ of $N$, and hence $\H x$ is normal
in $N$.  (i) and (ii) are elementary consequences of the first part.
  \proofend


\section{A formula for the product and relations for the Hecke algebra}
  From now on we assume that $(G,H)$ is a Hecke pair such that $H$ is
subnormal in $G$.  One of our main goals is to obtain what we believe
are the cleanest formulas ever for the product
  $\s x\s y$
  of generators of $\HGH$.  See also \scite{\Krieg}{I.4.4} and
\cite{\QuiggKaliLand}.

We begin by studying certain aspects of double cosets in a little more
detail.  As usual let us view double cosets as the orbits of the
action $\beta$ of $H\times H$ on $G$ given by
  $$
  \beta_{(h,k)}(x)=hxk\inv
  \for (h,k)\in H\times H\for x\in G.
  $$

  Given double cosets $HxH$ and $HyH$, observe that their product
  $$
  HxHyH
  $$
  is invariant under $\beta$ and hence may be written as the disjoint
union of orbits of the form
  $
  HxhyH,
  $
  for certain elements $h$ in $H$.
  Observe moreover that, for $h,k\in H$, one has that
  $$
  HxkyH = HxhyH
  \eqno{(\dagger)}
  $$
  if and only if
  $$
  xky\in  HxhyH \ \Leftrightarrow \
  k\in x\inv HxhyHy\inv = \H xh\H{y\inv}.
  $$

  Since $H$ is subnormal  in $G$ we have that  $\H xh=h\H x$ and
hence  $(\dagger)$ holds if and only if $h$ and $k$ define
the same coset modulo $H\cap\H x\H{y\inv}$.
  We therefore have:

  \state Proposition 
  \relabel\DoubleInDoubleProp
  Suppose that $H$ is a subnormal subgroup of a group $G$.
  Given $x,y\in G$, let $S_{x,y}$ be a family of representatives for
  $\rq H{(H\cap \H x\H{y\inv})}$.  Then
  $$
  HxHyH = \dot{\bigcup_{h\in S_{x,y}}} \ HxhyH.
  \eqmark\DoubleInDouble
  $$
  Conversely, if $S_{x,y}$ is any subset of $H$ such that
\lcite{\DoubleInDouble} holds then it is a family of representatives
for $\rq H{(H\cap \H x\H{y\inv})}$.

  \proof
  If $h,k\in H$ notice that
  $$
  \def\bi{\ \Leftrightarrow\ }
  HxhyH =   HxkyH \bi
  xhy \in   HxkyH \bi
  h \in   x\inv HxkyHy\inv = \H xk\H{y\inv}.
  $$
  Under the hypothesis that $H$ is subnormal we have that $\H xk=k\H
x$ so the above holds if and only if $hk\inv\in H\cap \H x\H{y\inv}$.
  \proofend

We are now ready to prove an important result, namely
that the product of two elements $\s x$ and $\s y$, corresponding to
the double cosets $HxH$ and $HyH$, is the average of the $\s z$ for
the double cosets $HzH$ which make up $HxHyH$.

  \state Theorem 
  \relabel\MainTheoremOnHeckeProduct
  Let $(G,H)$ be a Hecke pair such that $H$ is subnormal in $G$.
Given $x,y\in G$, let $S_{x,y}$ be any subset of $H$ such that $HxHyH
= \dot{\bigcup}_{h\in S_{x,y}} HxhyH$.  Then
  $$
  \s x\s y =
  {1\over |S_{x,y}|}\sum_{h\in S_{x,y}}\s{xhy}.
  $$

  \proof Consider the following diagram in which all horizontal maps
are defined to be the inclusion of the group appearing in the
corresponding numerator, moded out by the corresponding denominators:
  $$
  \def\rght{\ \longrightarrow \ }
  \matrix{
  0 & \rght &
  \ds {H\cap \H x \H{y\inv} \over H\cap\H x} & \rght &
  \ds {H \over H\cap\H x} & \rght &
  \ds {H \over H\cap\H x\H{y\inv}} & \rght & 
  0 \cr\cr\cr
  && \kern 25pt\Big\downarrow {Ad}_{y\inv} \cr\cr\cr
  0 & \rght &
  \ds{\H y \cap \H{xy}H\over \H y\cap\H{xy}} & \rght &
  \ds{H\H{xy} \over \H{xy}} & \rght &
  \ds{H\H{xy}\H y \over \H{xy}\H y} & \rght &
  0 \cr\cr\cr
  &&&&&&\Big\|\cr\cr\cr
  0 & \rght &
  \ds{H\H y \cap \H{xy}\H y\over \H y} & \rght &
  \ds{H\H y \over \H y} & \rght &
  \ds{H\H{xy}\H y \over \H{xy}\H y} & \rght &
  0 
  }
  $$
  It is elementary to check that all rows are exact and all vertical
maps are isomorphisms.

  We will refer to these groups by the cardinal points so that for
instance ${H\cap \H x \H{y\inv} \over H\cap\H x}$ will be called the
northwest group.  

  \def\NW{A}
  \def\NE{S_{x,y}}
  \def\SW{C}
  \def\SE{D}

  Recalling that $\NE$ is a family of representatives for the
northeast group by \lcite{\DoubleInDoubleProp}, let $\NW$ be a family
of representatives for the northwest group, so that the set
  $$
  S_x := \NE \NW = \{ba: b\in \NE,\ a\in \NW\}
  $$
  is a family of representatives for the north group.  It is also
worth noticing that for distinct pairs $(a_1,b_1)$ and $(a_2,b_2)$ in
$\NW\times\NE$ one has that $b_1a_1\neq b_2a_2$, so that $|S_x|=|\NW| |\NE|$.

  Similarly let $\SW$ and $\SE$ be families of representatives for the
southwest and southeast groups, respectively.
Notice that the equivalence class of an element in $\SE$ is unaltered
upon multiplication by an element from $\H{xy}\H y$ so we may supose
that 
  $$
  \SE\subseteq H.
  \eqmark\SEH
  $$
  It follows that $\SW\SE\subseteq H\H y$ and it
is then clear that $T_y:=\SW\SE$ is a family of representatives for
the south group.  As above $|T_y|=|\SW|\ |\SE|$.

Also, observe that the equivalence class of an element in $\SW$ is
unaffected under multiplication by an element from $\H y$ and hence we may
assume that
  $$
  \SW \subseteq  \H{xy}.
  \eqmark\SWHxy
  $$

  Using  \lcite{\AlternateHeckeOperator} for the description of $\s y$
we have
  $$
  \s x\s y\calcat{\dH{}} =
  {1\over |S_x||T_y|}\sum_{h\in S_x}\sum_{k\in T_y} \dH{xhyk}=
  {1\over |\NW||\NE||\SW||\SE|}\sum_{b\in \NE}\sum_{a\in
\NW}\sum_{c\in \SW}\sum_{d\in \SE} \dH{xbaycd}.
  $$
  For each $b$, $a$, $c$, and $d$ as above notice that
  $$
  xbayc d =
  (xbay)c(xbay)\inv xbay d =
  c'xbayd,
  $$
  where 
  $$
  c' =
  (xbay)c (xbay)\inv {\buildrel (\SWHxy) \over \in}
  xbay \H{xy} y\inv a\inv b\inv x\inv = H.
  $$
  Therefore 
  $
  Hxbaycd = Hxbayd,
  $
  so that
  $$
  \s x\s y\calcat{\dH{}} =
  {1\over |\NW||\NE||\SE|}\sum_{b\in \NE}\sum_{a\in \NW}\sum_{d\in
\SE} \dH{xbayd}.
  $$
  Next write 
  $$
  Hxbayd =
  Hxbyy\inv ayd =
  Hxbya'd,
  $$
  where $a'= y\inv ay$.  Denoting by $\NW'=y\inv\NW y$, we see that $\NW'$
is a family of representatives for the west group and thus, by
\lcite{\SEH}, we have that 
  $$
  T_{xy}:=\NW'\SE
  $$ 
  is a subset of $H\H{xy}$ as well as a family of representatives for
the center group.  So
  $$
  \s x\s y\calcat{\dH{}} =
  {1\over |\NW'||\NE||\SE|}\sum_{b\in \NE}\sum_{a'\in \NW'}\sum_{d\in
\SE} \dH{xbya'd} \$=
  {1\over |\NE||T_{xy}|}\sum_{b\in \NE}\sum_{\ell\in T_{xy}}
\dH{xby\ell} =
  {1\over |\NE|}\sum_{b\in \NE}\s{xby}\calcat{\dH{}}.
  \proofend
  $$

Based on examples we have been able to determine that
the above product formulas do not hold for general Hecke pairs.  We
therefore leave open the following:

 \sysstate {Question}{\rm}{\relabel\OpenQuestion For which Hecke pairs do the
product formulas of \lcite{\MainTheoremOnHeckeProduct} hold?  Do they
hold when $H$ is protonormal?}

Back to the subnormal realm we obtain the following universal property of Hecke algebras:

  \state Theorem
  \label\UniversalProperty
  Suppose that  $(G,H)$ is a Hecke pair such that $H$ is subnormal in
$G$ and  let $\letterfort$ be  any map from $G$ into  a
unital $\Fd$-algebra $B$ such that $\t 1=1$, and for every $x,y\in G$ and for every
finite set $S_{x,y}\subseteq G$ such that
  $$
  HxHyH = \dot{\bigcup_{h\in S_{x,y}}} \ HxhyH,
  $$
  one has that
  $$
  \t x\t y =
  {1\over |S_{x,y}|}\sum_{h\in S_{x,y}}\t{xhy}.
  \eqmark\GeneralLinearCombination
  $$
  Then:
  \izitem 
  \zitem There exists a unique unital homomorphism $\phi:\HGH\to B$
such that $\phi(\s x)=\t x$, for all $x\in G$.
  \zitem If moreover $B$ is a *-algebra and $\t{x\inv}=\t x^*$, for all
$x$ in $G$, then $\phi$ is a *-homomorphism.

  \proof 
  Given $x$ in $G$ and $h$ in $H$ we have that
  $$
  \matrix{
  HxHhH = HxH & = &  Hxh\inv hH\cr
  & = & Hx1hH, \hfill
  }
  $$
  which means that the singletons
$\{h\inv\}$ and $\{1\}$ are acceptable choices for $S_{x,h}$.  Therefore we have by
\lcite{\GeneralLinearCombination} that
  $$
  \matrix{
  \t x\t h  & = &  \t{xh\inv h}\cr
  & = & \t{x1h},\hfill
  }
  $$
  which implies that $\t x=\t{xh}$.  Beginning with $HhHxH=HxH$ one
may similarly conclude that $\t x=\t{hx}$.  It therefore follows that
$\letterfort$ is a doubly invariant function on $G$.  Employing
\lcite{\LinearBasis} we therefore see that there exists a unique \stress{linear map}
$\phi:\HGH\to B$ such that $\phi(\s x)=\t x$, for all $x\in G$.

So  we need only prove that $\phi$ is a homomorphism in order to
establish (i).  In order to do this it is obviously enough to prove
that $\phi(\s x\s y)=\phi(\s x)\phi(\s y)$, for all $x$ and $y$ in $G$.
Given $S_{x,y}$ as in the statement  we have
  $$
  \phi(\s x\s y) \={(\MainTheoremOnHeckeProduct)}
  {1\over |S_{x,y}|}\sum_{h\in S_{x,y}}\phi(\s{xhy}) =
  {1\over |S_{x,y}|}\sum_{h\in S_{x,y}}\t{xhy} =
  \t x \t y =
  \phi(\s x)  \phi(\s y),
  $$
  proving our claim that $\phi$ is a homomorphism.  

  If $B$ is a *-algebra and $a\in\HGH$ is the finite sum
  $a = \sum_{x\in G} \lambda_x \s x$, then
  $$
  \phi(a^*) =
  \phi\Big(\sum_{x\in G} \conj{\lambda_x} \s x^* \Big) \={(\StarAlgebraStructure)}
  \phi\Big(\sum_{x\in G} \conj{\lambda_x} \si x \Big) =
  \sum_{x\in G} \conj{\lambda_x} \t{ x\inv}  \$=
  \sum_{x\in G} \conj{\lambda_x} \t x^* =
  \Big(  \sum_{x\in G} \lambda_x \t x\Big)^* =
  \phi(a)^*.
  \proofend
  $$

Putting together \lcite{\MainTheoremOnHeckeProduct} and
\lcite{\UniversalProperty} we arrive at the following presentation of
the Hecke algebra.

\state Theorem  
  \relabel\Relations
  Let $(G,H)$ be a Hecke pair such that $H$ is subnormal in $G$.  Then
the Hecke algebra $\HGH$ admits the following presentation in the
category of unital $\Fd$-algebras: 
  \iItem 
  \Item{\rm (a)} GENERATORS:  any set indexed by  $G$,
say $\{\t x:x\in G\}$,
  \Item{\rm (b)} RELATIONS: 
  \def\zitemitem {\global \advance \zitemno by 1 \smallskip \itemitem{{\rm(\romannumeral
\zitemno)}}}
  \izitem
  \zitemitem $\t 1 =1$,
  \zitemitem $\t x\t y = \ds {1\over |S_{x,y}|}\sum_{h\in
S_{x,y}}\t{xhy}$, whenever $S_{x,y}$ is a subset of $H$ such that $HxHyH
= \dot{\bigcup}_{h\in S_{x,y}} HxhyH$.
  \medskip\noindent
  If we add
  \zitemitem $\t{x\inv} = \t x^*$, for every $x$ in $G$,
  \medskip\noindent
  we arrive at a presentation of\/ $\HGH$ in the category of unital
*-algebras over $\Fd$.

Let us now study some simple properties shared by maps $\letterfort$
satisfying the above relations:

\state Proposition
  \label\FurtherProperties
  Let $(G,H)$  be a Hecke pair with $H$ subnormal in $G$ and let  $B$
be a unital $\Fd$-algebra.   Given any  map $\letterfort:G\to B$ 
satisfying \lcite{\Relations.b.i-ii} we have
  \izitem
  \zitem $\letterfort$ is a partial representation,
  \zitem if $xH\subseteq Hx$, then $\t x\t{x\inv}=1$, and $\t x \t y =
\t{xy}$, for all $y$ in $G$,
  \zitem if $Hx\subseteq xH$, then $\t{x\inv}\t x=1$, and $\t y \t x =
\t{yx}$, for all $y$ in $G$,
  \zitem if $x$ lies in the normalizer of $H$ 
then $\t x$ is invertible, $(\t x)\inv=\t{x\inv}$, $\t x\t y=\t{xy}$
and $\t y\t x=\t{yx}$, for all $y$ in $G$,
  \zitem $\t x = \t{hxk}$, for all $h$ and $k$ in $H$.

  \proof
  In order to prove (i) let $\phi:\HGH\to B$ be the homomorphism given
by \lcite{\UniversalProperty.i}.  Then for every $x$ and $y$ in $G$ we have
  $$
  \t{x\inv}\t x\t y = 
  \phi(\s{x\inv})\phi(\s x)\phi(\s y) = 
  \phi(\s{x\inv}\s x\s y) \= {(\MainPrep)}
  \phi(\s{x\inv}\s{xy}) = 
  \t{x\inv}\t{xy},
  $$
  while a similar argument proves that 
  $\t x\t y\t {y\inv} = \t {xy}\t {y\inv}$.

  Supposing that $xH\subseteq Hx$, we have that $HxHx\inv H=H=H1H$, so we
may take $S_{x,x\inv}=\{1\}$ in \lcite{\Relations.b.ii} to conclude
that $\t x\t{x\inv}=\t 1 = 1$.  
  Moreover
  $$
  \t x \t y =   \t x \t{x\inv} \t x \t y = 
  \t x \t{x\inv} \t {xy} = 
  \t {xy}.
  $$
  Clearly (iii) follows from (ii) by taking inverses, while (iv) follows
from  (ii) and (iii).
  As for (v) first notice that taking $x\in H$ and $y=1$ in
\lcite{\Relations.b.ii} we have that
  $HxHyH = H = Hxx\inv yH$ so, taking  $S_{x,y}=\{x\inv\}$ we get
  $$
  \t x =  \t x \t 1 = \t{xx\inv 1} = 1,
  $$
  proving that $\letterfort$ is constantly equal to 1 on $H$.  Since $H$
obviously normalizes itself we have that (v) follows from (iv).
  \proofend


  \section{Hecke algebra as a crossed product}
  Throughout this section we fix a Hecke pair $(G,H)$ and a subgroup $N$
of $G$ such that $H\normal N\normal G$, in which case $H$ is
necessarily subnormal
in $G$.  Our goal will be to show that there exists a twisted partial
action of $\rq GN$\fn{We seem to be irremediably biased towards right
coset spaces so we will keep using the notation for right cosets even
when they coincide with left cosets.} on the group algebra $\F{\rq
NH}$ such that the corresponding crossed-product is isomorphic to
$\HGH$.

By \lcite{\FurtherProperties.iv} we have that the restriction of
$\letterfors$ to $N$ is a global (as opposed to partial)
representation of $N$ on $\HGH$.  Since $\letterfors$ vanishes on $H$
we in fact get a group representation of $\rq NH$ on $\HGH$ which maps
each right (= left) coset $Hn$ in $\rq NH$ to $\s n$.

\state Proposition
  The homomorphism $\iota:\F{\rq NH}\to\HGH$ obtained by linearizing
the above representation of $\rq NH$ is injective.

  \proof
  Given $n\in N$ it is evident that $\{1\}$ is a family of
representatives for $\rq H{(\HH n)}$, so
  $
  \s n\calcat{\dH{}} =
  \dH{n}.
  $
  We thus see that for a general element
  $$
  a = \sum_{i=1}^n \lambda_i \dH{n_i} \in \F{\rq NH}
  $$
  (by abuse of language  we denote by $\dH n$ the canonic basis
elements of $\F{\rq NH}$ as well)
  one has that
  $$
  \iota(a) =   
  \sum_{i=1}^n \lambda_i \iota(\dH{n_i}) =
  \sum_{i=1}^n \lambda_i \s {n_i},
  $$
  whence
  $$
  \iota(a)\Calcat{\dH{}} =   
  \sum_{i=1}^n \lambda_i \s {n_i} \calcat{\dH{}} =   
  \sum_{i=1}^n \lambda_i \dH{n_i},
  $$
  from which the statement follows.
  \proofend

Using $\iota$ we will identify, from now on, $\F{\rq NH}$ with a
sub-algebra of $\HGH$, namely the linear span of the set $\{\s n:
n\in N\}$.

  \state Lemma
  \label\InFNH
  Let $x,y\in G$ be such that $xy\in N$.  Then $\s x\s y\in
\F{\rq NH}$.
  
  \proof
  By \lcite{\MainTheoremOnHeckeProduct} it is enough to show that
$xhy\in N$, for every $h\in H$.
  Let $n=xy$, so that $y=x\inv n$. Thus, given $h\in H$,
we have
  $$
  xhy = xhx\inv n \ \in\  xHx\inv n \ \subseteq \ xNx\inv n = N n =N.
  \proofend
  $$
  
  \state Lemma
  \label\CentralIdempotent
  For every $x\in G$ one has that $\e x:= \s x\s {x\inv}$ is a central
idempotent in $\F{\rq NH}$.  Moreover if $Nx=Ny$ then $\e x=\e y$.

  \proof From \lcite{\MainPrep} it follows that $\e x$ is an
idempotent and  from \lcite{\InFNH}, that $\e x\in \F{\rq NH}$.

Let $S_x$ be a family of representatives for $\rq{H}{(\HH x)}$.
Plugging $y=x\inv$ in  \lcite{\DoubleInDoubleProp}
we have that   $HxHx\inv H = \dot{\bigcup}_{h\in S_{x}} \ Hxhx\inv H$
so, by \lcite{\MainTheoremOnHeckeProduct},
  $$
  \e x =\s x\s {x\inv} =
  {1\over |S_{x}|}\sum_{h\in S_{x}}\s{xhx\inv}.
  $$
  In order to prove that $\e x$ is central it is enough to show that
$\e x$ commutes with $\s n$, for all $n\in N$.  For this observe that
$m:= x\inv n x \in N$, so $Ad_m$ is an inner automorphism of $N$,
which therefore leaves invariant the normal subgroups $H$ and $\H x$.
We conclude that $Ad_m(S_x)$ is another family of representatives
for $\rq{H}{(\HH x)}$ so we can alternatively compute $\e x$ as
  $$
  \e x =\s x\s {x\inv} =
  {1\over |S_{x}|}\sum_{h\in S_{x}}\s{xmhm\inv x\inv} \$=
  {1\over |S_{x}|}\sum_{h\in S_{x}}\s{n xhx\inv n\inv}
\={(\FurtherProperties.iv)}
  {1\over |S_{x}|}\sum_{h\in S_{x}}\s n\s{xhx\inv}\s{n\inv} =
  \s n \e x\s{n\inv} =
  \s n \e x(\s n)\inv,
  $$ 
  proving that $\e x$ commutes with $\s n$.
  If $Nx=Ny$ we may write $y=nx$, with $n\in N$, so
  $$
  \e y =
  \e{nx} =
  \s{nx}\s{x\inv n\inv} \={(\FurtherProperties.iv)}
  \s n\s x\s{x\inv}\s{n\inv} = 
  \s n\e x\s{n\inv} = 
  \e x.
  \proofend 
  $$

  \definition For each $x\in G$ we will let:
  \izitem 
  \zitem $D^x$ be the ideal of $\F{\rq NH}$ generated by $\e x$, that
is  $D^x = \e x\F{\rq NH}$,
  \zitem $\psi_x$ be the linear operator on $\HGH$ given by
  $$
  \psi_x : a \in \HGH \mapsto \s x a\s{x\inv} \in \HGH.
  $$

By the last part of \lcite{\CentralIdempotent} it is clear that $D^x$
only depends on the class of $x$ in $\rq GN$.  If $t\in \rq GN$ we
will therefore denote by
  $$
  D_t := D^x,
  \eqmark\DefDGamma
  $$
  where $x$ is any element of $G$ such that $Nx=t$, so 
$D_t$ is independent of the choice of $x$.

  \state Proposition
  \label\PartialAutos
  For every  $x\in G$ one has that $\psi_x(\F{\rq NH})= D^x$.
Moreover the restriction of $\psi_x$ to $D^{x\inv}$ is an isomorphism
onto $D^x$.
  
  \proof In order to verify that $\psi_x(\F{\rq NH})\subseteq D^x$ it is enough to show
that $a:= \s x \s n\s{x\inv} \in D^x$, for all $n\in N$.  Notice that
  $$
  a =
  \s x\s n\s{x\inv} \={(\FurtherProperties.iv)}
  \s {x}\s{nx\inv}\ {\buildrel (\InFNH) \over \in}\ \F{\rq NH}.
  $$
  Since $a=\e x a$, by \lcite{\MainPrep}, we conclude that
$a\in D^x$.

Observe that for $a\in D^{x\inv}$ we have
  $$
  \psi_{x\inv} (\psi_x(a)) =
  \s{x\inv} \s x a \s{x\inv}\s x =
  \e{x\inv} a  \e{x\inv}  = a
  $$
  from which it follows that $\psi_x$ is a bijection from $D^{x\inv}$
to $D^x$.  From this we also obtain  that  $\psi_x(\F{\rq NH})=D^x$.
  Finally, in order to show that the restriction of $\psi_x$
to $D^{x\inv}$ is multiplicative, let $a,b\in D^{x\inv}$.  Then
  $$
  \psi_x(ab) =
  \s x ab \s{x\inv} =
  \s x a\e{x\inv}b \s{x\inv} =
  \s x a\s{x\inv}\s xb \s{x\inv} =
  \psi_x(a)\psi_x(b).
  \proofend
  $$

Fix, once and for all, a section $\letterforr$ for the quotient map
$\pi:G\to \rq GN$, that is, $\letterforr$ is a map (not necessarily a
homomorphism) from $\rq GN$ to $G$ such that $\pi\circ \letterforr$ is
the identity map on $\rq GN$.  For the special case of the coset $N1$
we will force the choice
  $$
  \r {N1} = 1.
  $$

  Given $r,s\in \rq GN$, observe that 
  $$
  \pi\big(\r r \r s \r{rs}\inv\big) = rs(rs)\inv = 1,
  $$
  so the element 
  $
  \r r \r s \r{rs}\inv
  $
  lies in $N$.

  \definition
  \label\DefineCocycle
  For every $r$ and $s$ in $\rq GN$ we let 
  $$
  \w rs = \s{\r r \r s \r{rs}\inv}.
  $$
  
  Clearly $\w rs$ is an invertible element in $\F{\rq NH}$ by
\lcite{\FurtherProperties.iv}.

  \state Lemma
  \label\LemmaOnCocycles
  Given  $r$ and $s$ in $\rq GN$,  let  $x=\r r$, $y=\r s$, and  $z = \r{rs}$. Then
  \izitem    
  \zitem $ \s x \s y \ei y =  \w rs \s z \ei y$,
  \zitem $ \ei y \si y \si x = \ei y  \si z (\w rs)\inv.$

  \proof Letting $n=  \r r \r s \r{rs}\inv = xyz\inv$ we have that
$xy=nz$ and $\w rs = \s n$.  So
  $$
  \s x \s y \ei y =
  \s x \s y \si y \s y =
  \s{xy}\si y\s y =
  \s{nz}\ei y \={(\FurtherProperties.iv)}
  \s n\s z\ei y =
  \w rs \s z\ei y.
  $$  
  As for (ii) we have
  $$
  \ei y \si y \si x =
  \si y \s y \si y \si x =
  \si y \s y \s{y\inv x\inv} \$=
  \ei y  \s{z\inv n\inv} \={(\FurtherProperties.iv)}
  \ei y  \si z \si n \={(\FurtherProperties.iv)}
  \ei y  \si z (\s n)\inv  =
  \ei y  \si z (\w rs)\inv.
  \proofend
  $$

\state Theorem For each $t\in\rq GN$, let $D_t$ be as in
\lcite{\DefDGamma} and let 
$\theta_t$ be the isomorphism from $D_{t\inv}$ to $D_t$
given by restricting $\psi_{\r{t}}$ to $D_{t\inv}$ as in
\lcite{\PartialAutos}.  Then the triple
  $$
  \Big(
    \{D_t\}_{t\in\rq GN}, \quad
    \{\theta_t\}_{t\in\rq GN}, \quad
    \{\w ts\}_{t,s\in\rq GN}
  \Big)
  $$
  is a twisted partial action of $\rq GN$ on $\F{\rq NH}$.
  
  \proof
  During the course of this prove we will let $A:=\F{\rq NH}$.

  Since $\s 1=1$,  it is evident that $D_1=A$ and $\theta_1$ is the
identity map on $A$. 
  In order to verify \lcite{\DefineTPA.ii} let $r,s\in\rq GN$, and put
$x=\r r$ and $y=\r s$.  So
  $$
  \theta_r(D_{r^{-1}}\cap D_s) =
  \s x(\e {x\inv}A \cap \e y A) \s{x\inv} =
  \s x(\e {x\inv} \e y A) \s{x\inv} =
  \s x\si x\s x \s y\si y A \s{x\inv} \$=
  \s {xy}\si y A \s{x\inv} =
  \s {xy}\si y A \ei x\si x =
  \s {xy}\si y \si x\s x A\si x =
  \s {xy}\s{y\inv x\inv}\s x A\si x \$=
  \e {xy}\s x A\si x \={(\PartialAutos)}
  \e {xy}D^x =
  D^{xy}\cap D^x =
  D_{rs}\cap D_r.
  $$

  As for \lcite{\DefineTPA.iii} let $x=\r r$ and $y=\r s$ as above and
put $z = \r{rs}$.  Take 
$a\in D_{s^{-1}} \cap D_{s^{-1}r^{-1}}$, which we may clearly suppose
has the form
  $$
  a = \ei y \e{(xy)\inv}\s n,
  $$
  where $n\in N$.  Then
  $$
  \def\.{\hfill}
  \def\eq{= \vrule height 12pt width 0pt}
  \matrix{\theta_r(\theta_s(a)) 
            & \eq & \s x \s y \kern 1.6cm a\kern 1.6cm \si y \si x \cr
            & \eq &  \s x \s y \. \ei y \e{(xy)\inv}\s n \.\si y \si x \cr
            & \={(\LemmaOnCocycles.i)} &
                   \.\w rs\s z \. \ei y \e{(xy)\inv}\s n \si y \si x \cr
            & \eq &  \w rs\s z \. \e{(xy)\inv} \ei y\. \s n \si y \si x \cr
            & \eq &  \w rs\s z \. \ei y \s n \si y \si x \cr
            & \eq &  \w rs\s z \. \s n \ei y \.\si y \si x \cr
            & \={(\LemmaOnCocycles.ii)} &
                   \w rs\s z \s n \. \ei y \si z(\w rs)\inv\. \cr
            & \eq &  \w rs\s z \. \e{(xy)\inv}\. \s n \ei y \si z(\w rs)\inv \cr
            & \eq &  \w rs\s z \.a\. \si z(\w rs)\inv \cr
            & \eq &  \w rs\.\theta_{rs}(a)\.(\w rs)\inv.}
  $$

  The forced choice of $\r {N1} = 1$ clearly gives
\lcite{\DefineTPA.iv} so it remains to check \lcite{\DefineTPA.v}.  So
let $r,s,t\in\rq GN$ and $a \in D_{r^{-1}}\cap D_s\cap D_{st}$.
Put  $x=\r r$, $y=\r s$, $z=\r t$, 
$\alpha=\r{rs}$, 
$\beta=\r{st}$, and $\gamma=\r{rst}$.
  We then have
  $$
  a = \ei x a = \e y a = \e \beta a,
  $$
  while
  $$
  \w rs= \s{x y \alpha\inv}, \quad
  \w st = \s{y z \beta\inv}, \quad
  \w r{st}= \s{x \beta \gamma\inv}, \quad
  \w {rs}t= \s{\alpha z \gamma\inv}.
  $$
  Therefore we have
  $$
  \theta_r(a\w st)\w r{st}=
  \s xa  \s{y z \beta\inv}  \si x  \s{x \beta \gamma\inv} \={(!)}
  \s xa  \s{(y z \beta\inv) x\inv  (x \beta \gamma\inv)} \$= 
  \s xa  \s{y z \gamma\inv} =
  \s x a \s {x\inv (x y \alpha\inv)(\alpha z \gamma\inv)} \={(!)}
  \s x a \si x \s{x y \alpha\inv}\s{\alpha z \gamma\inv} =
  \theta_r(a)\w rs\w{rs}t.
  $$
  Observe that the passages marked ``$(!)$'' are justified by
\lcite{\FurtherProperties.iv} and the fact that the elements $y z \beta\inv$, $x \beta
\gamma\inv$, $xy\alpha\inv$, and $\alpha z\gamma\inv$ lie in $N$.
  \proofend

  \state Theorem 
  \relabel\MainTheoremOnCrossedProduct
  The crossed product
  $$
  \F{\rq NH}\crossproduct \rq GN
  $$
  relative to the above twisted partial action is isomorphic to the
Hecke algebra $\HGH$.

  \def\dcp#1{\delta_{#1}} 
  \proof Let 
  $$
  \Phi:   \F{\rq NH}\crossproduct \rq GN \to \HGH
  $$
  be the unique linear map such that 
  $$
  \Phi(a\dcp t) = a\s{\r t}
  \for t\in \rq GN\for a\in D_t.
  $$
  In order to show that $\Phi$ is multiplicative 
let $r,s\in \rq GN$ and take $a\in D_r$ and $b\in D_s$.  Putting $x=\r
r$ and $y=\r s$ we have
  $$
  \Phi(a\dcp r) \Phi(b\dcp s) =
  a\s x b\s y =
  \e xa\s x b\s y =
  \s x\si xa\s x b\s y =
  \s x\theta_r\inv(a)b\s y \$=
  \s x\theta_r\inv(a)b\ei x\s y =
  \s x\theta_r\inv(a)b\si x\s x\s y =
  \theta_r(\theta_r\inv(a)b)\s x\s y \$=
  \theta_r(\theta_r\inv(a)b)\s x\s y \ei y  = \ldots
  $$
  Putting $z = \r{rs}$ and applying \lcite{\LemmaOnCocycles.i} we find
that the above equals
  $$
  \ldots = 
  \theta_r(\theta_r\inv(a)b)\w rs \s z \ei y
\={(\BasicCommutingRelation)}
  \theta_r(\theta_r\inv(a)b)\w rs \e{zy\inv} \s z =\ldots 
  $$
  Notice that $\pi(zy\inv) = \pi(z)\pi(y)\inv = (rs)s\inv = r = \pi(x)$ which
implies that $Nzy\inv =Nx$.  Hence $\e{zy\inv} = \e
x$ by \lcite{\CentralIdempotent} and the above equals
  $$
  \ldots =
  \theta_r(\theta_r\inv(a)b)\w rs \e x \s z =
  \theta_r(\theta_r\inv(a)b)\e x \w rs \s z =
  \theta_r(\theta_r\inv(a)b) \w rs \s z.
  $$
  On the other hand, since 
  $
  (a\dcp r)(b\dcp s) =  
  \theta_r(\theta_r\inv(a)b) \w rs \dcp {rs},
  $
  we have that
  $$
  \Phi\big(  (a\dcp r)(b\dcp s) \big) =
  \theta_r(\theta_r\inv(a)b) \w rs \s z,
  $$
  proving that $\Phi$ is a homomorphism.  
  In order to prove that $\Phi$ is bijective we will now provide an
inverse for it based on the universal property
\lcite{\UniversalProperty} of the Hecke
algebra.

Consider the map 
  $$
  \tau:G\to\F{\rq NH}\crossproduct \rq GN
  $$
  given by 
  $$
  \tau(x) = \e x\s{x\r{\pi(x)}\inv} \dcp {\pi(x)}.
  $$
  In order to simplify the above expression we will often write it as
  $$
  \tau(x) = \e x\s n \dcp r
  $$
  where 
  \def\xt{\tilde x} \def\xti{\xt\inv}%
  \def\yt{\tilde y}%
  \def\zt{\tilde z}%
  $r=\pi(x)$, $\xt =\r{r}$, and $n=x\xt\inv$.  Observe that $n$ is
necessarily in $N$.  

\bigskip\noindent CLAIM: Given  $x$ and $y$  in $G$ let
  $$
  \matrix{
  r=\pi(x),\hfill & \xt =\r{r},\hfill & n=x\xt\inv,\hfill \cr
  s=\pi(y),\hfill & \yt =\r{s},\hfill & m=y\yt\inv,\hfill\cr
  &\zt=\r{rs},
  }
  \eqmark\Conventions
  $$
  Then
  $$
  \t x \t y =
  \s x \s y \s{\zt\inv} \dcp {rs}
  \eqmark\ClaimProdTau
  $$

  In fact we have 
  $$
  \t x \t y =
  \big(\e x\s n \dcp r\big) \big(\e y\s m \dcp s\big) =
  \theta_ r\big(\theta_ r\inv(\e x\s n) \e y\s m\big) \w rs\dcp {rs} \$=
  \s {\xt}\big(\s{\xti}(\e x\s n)\s\xt \e y\s m \big) \s{\xti} \w rs
\dcp {rs} = 
  \e{\xt}\e x\s n\s\xt \e y\s m\s{\xti} \s{\xt\yt\zt\inv} \dcp {rs} \$=
  \e x \s {n\xt} \e y\s m\s{\yt\zt\inv} \dcp {rs} =
  \e x \s {n\xt} \e y\s {m\yt\zt\inv} \dcp {rs} \$=
  \e x \s x \e y\s {y\zt\inv} \dcp {rs} =
  \s x \e y\s {y\zt\inv} \dcp {rs} =
  \s x \s y \si y\s {y\zt\inv} \dcp {rs} =
  \s x \s y \s{\zt\inv} \dcp {rs},
  $$
  proving our claim.  
  Next let us show that $\letterfort$ satisfies
\lcite{\GeneralLinearCombination}.  For this let $S_{x,y}$ be a family
of representatives for $\rq H{(H\cap \H x\H{y\inv})}$.  Using our
claim and \lcite{\MainTheoremOnHeckeProduct} we conclude that
  $$
  \t x \t y =
  {1\over |S_{x,y}|} \sum_{h\in S_{x,y}} \s {xhy} \s{\zt\inv} \dcp {rs}.
  $$
  On the other hand, in order to compute the right-hand-side of
\lcite{\GeneralLinearCombination}, namely the sum
  $$
  {1\over |S_{x,y}|}\sum_{h\in S_{x,y}}\t{xhy},
  $$
  we observe that $\pi(xhy) = \pi(x)\pi(y)=rs$, so that
$\r{\pi(xhy)}=\r{rs}=\zt$.  This implies that
  $$
  \t{xhy} =
  \e{xhy}\s{xhy\zt\inv} \dcp{rs} \$=
  \s{xhy}\s{(xhy)\inv}\s{xhy\zt\inv} \dcp{rs} =
  \s{xhy}\s{\zt\inv} \dcp{rs}.
  $$

  This shows that \lcite{\GeneralLinearCombination} holds and hence
by the universal property of $\HGH$ we
conclude that there exists a homomorphism
  $$
  \Psi:   \HGH \to \F{\rq NH}\crossproduct \rq GN 
  $$
  such that $\Psi(\s x) = \t x$, for all $x$ in $G$.
  We claim that $\Psi$ is the inverse of $\Phi$.  In fact, using
\lcite{\Conventions}, we have
  $$
  \Phi(\Psi(\s x)) =
  \Phi(\t x) =
  \Phi(\e x\s n \dcp r) =
  \e {x} \s {n}\s{\r r} =
  \e {x} \s {n}\s \xt \= {(\FurtherProperties.iv)}
  \e {x} \s {n\xt} = 
  \e {x} \s {x} = 
  \s x.
  $$
  This shows that $\Phi\circ\Psi$ is the identity on $\HGH$.
To show that  $\Psi\circ\Phi$  is also the identity on
  $\F{\rq NH}\crossproduct \rq GN$
  it is clearly enough to check that
  $\Psi(\Phi(a))=a$, for every a in $\F{\rq NH}\crossproduct \rq GN$ of
the form $a = \e x\s p\dcp r$, where  $p$ is in $N$, and
we are again using \lcite{\Conventions}.
  We have
  $$
  \Phi(\e x\s p\dcp r) =
  \e x\s p\s\xt =
  \s p\e x\s\xt =
  \s p\e \xt\s\xt =
  \s p\s\xt =
  \s {p\xt}.
  $$
  Thus
  $$
  \Psi(\Phi(\e x\s p\dcp r)) =
  \Psi(\s {p\xt}) =
  \t {p\xt} =
  \e {p\xt}\s{{p\xt}\r{\pi({p\xt})}\inv} \dcp {\pi({p\xt})} =
  \e x\s p \dcp r,
  $$
  concluding the proof.
  \proofend

As discussed in the introduction,
many authors have considered the problem of describing Hecke algebras
as crossed products by semigroups, assuming that $G$ has a semi-direct product
structure.  It is easy to see that $G$ can be written as a semi-direct
product $G=N\crossproduct K$, where $K$ is another group (necessarily
$K=\rq GN$), if and only
if there exists a section $\letterforr: \rq GN \to G$ for the quotient
map which is a group homomorphism.  In this case notice that the
cocycle $w$ defined in \lcite{\DefineCocycle} becomes trivial.
We therefore have:

  \state Corollary
  \relabel\UntwistedCrossProduct
 Let $G=N\crossproduct K$ be a semidirect product of groups and let
$H$ be a normal subgroup of $N$ such that $(G,H)$ is a Hecke pair\fn{See
\scite{\LaLaOne}{Proposition 1.7} for sufficient conditions for
$(G,H)$ to be a Hecke pair.}.  Then there is an (untwisted) partial
action of $K$ on $\F{\rq NH}$ such that
 $$
 \HGH \simeq   \F{\rq NH}\crossproduct K.
 $$


  \section{Hecke C*-algebras}
  In this section we take $F$ to be the field of complex numbers and
consider the existence of a maximum C*-norm on $\HGH$.  See the
introduction for references to similar results in the literature.  The
completion of $\HGH$ relative to this norm, when it exists, is
sometimes called the Hecke C*-algebra of the pair $(G,H)$ and its
*-representation theory is equivalent to the *-representation theory of
$\HGH$.  
Observe that by \lcite{\CompareStars} it does not matter which
involution we take on $\HGH$.
  
  \state Proposition Let $F=\C$ and let $(G,H)$ be a Hecke pair with $H$ protonormal
in $G$.  Then there exists a maximum C*-norm on $\HGH$.

  \proof
  Given $a\in\HGH$ let $\|a\|$ be  defined as the supremum of
$\|\pi(a)\|$, where $\pi$ ranges in the set of all *-representations of
$\HGH$.  To see that $\|a\|$ is finite
write $a$ as a finite sum
  $a = \sum_{x\in G} a_x\s x$.  Observe that if  $\pi$ is any *-representation of
$\HGH$ then, given that
  $$
  \s x \s x^* \s x = 
  \s x \si x \s x =
  \s x \s{x\inv x} =
  \s x \s 1 =
  \s x,
  $$
  we see that 
  $\pi(\s x)$ is a partial
isometry and hence $\|\pi(\s x)\|\leq 1$.  It follows that
  $$
  \|\pi(a)\|\leq 
  \sum_{x\in G} |a_x|\|\pi(\s x)\| \leq
  \sum_{x\in G} |a_x|.
  $$
  This proves that $\|a\|\leq   \sum_{x\in G} |a_x|$ and hence $\|a\|$
is finite as claimed. It is now easy to see that $\|\cdot\|$ defines a
C*-norm which dominates all others.
  \proofend

The completion of $\HGH$ relative to this norm is a C*-algebra
sometimes denoted by $C^*_u(G,H)$ and called called the \stress{full
Hecke C*-algebra} of the pair $(G,H)$.  It is elementary to see that
the *-representation theory of this algebra coincides with that of
$\HGH$.

On the other hand, as some authors have already done, one could consider
the reduced Hecke C*-algebra $C^*_r(G,H)$, namely the completion of
$\HGH$, normed as operators on the inner-product space defined in
\lcite{\NewInnerProduct}.  The question as to whether $C^*_u(G,H)$
coincides with $C^*_r(G,H)$ is then at least as rich as the
corresponding question for group C*-algebras.

  \section{A possible generalization of Hecke algebras}
  In this short section we wish to propose a generalization for the
definition of Hecke algebras for a group-subgroup pair $(G,H)$ which
is not a Hecke pair, namely, such that not all double cosets are
finite union of right cosets.

Initially observe that the relations \lcite{\Relations.b} make sense
as long as every ``triple coset'' $HxHyH$ is a finite union of double
cosets.  One could then be tempted to say that the pair $(G,H)$ is a
\stress{pseudo Hecke pair} if for every $x$ and $y$ in $G$ this
finiteness condition holds.

However observe that at least in the case of a subnormal $H\subseteq
G$, we have by \lcite{\DoubleInDoubleProp} that $HxHy H$ is a finite
union of double cosets if and only if $\rq H{(H\cap \H x\H{y\inv})}$
is finite.  If this is so for every $x$ and $y$ then, plugging
$y=x\inv$ we conclude that $\rq H{(H\cap \H x)}$ is finite and hence
$HxH$ is a finite union of right cosets by \lcite{\DoubleInRightProp}.
In other words every pseudo Hecke pair is a true Hecke pair.

However there is a lesson to be learned from \cite{\infinoa} which
could perhaps yield a true generalization.  That lesson is that, when
a collection of relations involves summations, some of which refuse to
converge, it is sensible to simple ignore the divergent ones.
A well known instance of this phenomenon takes place when one
considers  Cuntz algebras.  The relation 
  ``$\sum_{i=1}^nS_iS_i^* = 1$'' in the usual presentation of 
${\cal O}_n$ is simply ignored in the definition of 
${\cal O}_\infty$.

One could then risk the following:
  
  \definition 
  \relabel\GeneralizedHeckeAlgebra
  Let $H$ be a subnormal subgroup of a group $G$.  The generalized
Hecke algebra $\wt\HGH$ is the universal $F$-algebra generated by a
collection of elements $\{\s x: x\in G\}$ subject to the relations
declaring that $\letterfors$ is a partial representation in
addition to the following: whenever
  $HxHyH$ happens to be a finite union of double cosets (and only in
this case) we require that \lcite{\Relations.b.ii} holds as well.

  While we have nothing of interest to say at the moment about the
algebra so defined, it is not hard to give an example of a
group-subgroup pair $(G,H)$ which is not a Hecke pair although there
are many pairs of elements $x$ and $y$ for which $HxHyH$ is a finite
union of double cosets.  Consider for example
  $$
  G = \pmatrix{1 & \R\hfill \cr 0  & \R_+ } =
  \left\{\pmatrix{1 & b\cr 0 & a }\in {GL}_2({\bf R}): a,b\in\R,\ a>0\right\},
  $$
  with $H = \pmatrix{1 & \Z\cr 0 & 1 }$.  
  If 
  $x = \axb{b}{a}\in G$, it is easy to see that $\H x=\axb{a\Z}{1}$ hence, if $y
= \axb{d}{c}$ we have
  $$
  \H x\H{y\inv} = \axb{a\Z+d\inv\Z}{1}.
  $$
  Quite often one would have that $H\cap\H x\H{y\inv}=\{0\}$ in which
case there are infinitely many double cosets in $HxHyH$.  However if
the rational vector space generated by $a$ and $d\inv$ contains
a nonzero rational number then there will be an integral solution
$(n,m,p)$ to
the equation 
  $$ 
  an + d\inv m = p,
  $$
  with nonzero $p$, in which case $\axb{p\Z}{1}\subseteq H\cap \H
x\H{y\inv}$ so that $HxHyH$ will contain no more than $p$ double
cosets and then relation \lcite{\Relations.b.ii} would apply.

In Definition \lcite{\GeneralizedHeckeAlgebra} we have restricted
ourselves to the situation in which $H$ is subnormal in $G$ so that,
when $(G,H)$ is a Hecke pair, one recovers the usual Hecke algebra
$\HGH$ by \lcite{\Relations}.  However there does not seem to be any
immediate technical difficulty in adopting Definition
\lcite{\GeneralizedHeckeAlgebra} for a general (non-subnormal)
group-subgroup pair $(G,H)$ although this would most definitely depart
from the usual theory of Hecke algebras.

  \section{An example}
  In this section we shall give an example of a Hecke pair $(G,H)$,
such that $H$ is protonormal in $G$ but not subnormal.


Let $\P\subseteq \N$ be a set of prime numbers and let $\AP$ be the
subset of all rational numbers $n/m$, with $n,m\in \N$, $m\neq 0$,
such that no prime in $\P$ divides $m$.  It is clear that $\AP$ is a
subring of $\Rat$.  We will denote by $\AP^*$ the set of invertible
elements in $\AP$, so that a rational number $\xi$ lies in $\AP^*$ if
and only if $\xi=n/m$ and no prime in $\P$ divide either $n$ or $m$.

Denote by $G$ the group 
  $$
  G = \AXB{\Rat\hfill}{\Rat^*},
  $$
  meaning the set of all matrices  $\axb ba\in GL_2(\Rat)$, such that
$a\in\Rat^*=\Rat\setminus\{0\}$, and $b\in\Rat$,  and let $\HP$ be the subgroup
  $$
  \HP = \AXB{\AP}{\AP^*}.
  $$

  \state Proposition
  For any set $\P$ of primes one has that $(G,\HP)$ is a Hecke pair.
  
  \proof
  Let $x=\axb ba\in G$.  We shall prove that
$\rq{\HP}{(\HP\cap\HP^{x\inv})}$ is finite.  As a first step lets us
try to identify certain elements in  $\HP\cap\HP^{x\inv}$.
Given $h=\axb \eta\xi\in\HP$,
notice that $h\in\HP^{x\inv}$, if and only if $x\inv hx\in \HP$.  We have
  $$
  x\inv h x =
  \AXB {-ba\inv}{a\inv} \AXB \eta \xi \AXB ba =
  \AXB {\eta a+(1-\xi)b}{\xi}.
  \eqmark\ComputeConjugate
  $$
  Therefore $x\inv hx\in \HP$ if and only if 
  $$
  \eta a+(1-\xi)b\in\AP.
  \eqmark\LieInHPxinv
  $$

Since neither $a$ or $b$ have been  assumed to lie in $\AP$ their denominator
may contain factors in $\P$.  Factoring these out we may write
  $$
  a = {a_0\over p}
  \and
  b = {b_0\over q},
  $$
  where $a_0,b_0\in\AP$, and $p$ and $q$ are products of primes in
$\P$.  Writing $a = qa_0/pq$ and $b = pb_0/pq$, we may assume without
loss of generality that $p=q$.

  Let $\Z_{q}$ denote the ring $\Z/q\Z$.  Given  $\zeta\in\AP$ write
it in reduced form $\zeta=n/m$, so that  no prime in $\P$ divides $m$
and hence $\gcd(m,q)=1$ (greatest common divisor).  Therefore $m$ is
invertible modulo $q$ and hence it makes sense to set
  $$
  \phi(\zeta) = n m\inv (\mod\ q).
  $$
  This therefore gives a well defined  map $\phi:\AP\to\Z_q$, which can
be easily proven to be a homomorphism of rings.  Let $G_q$ be the
subgroup of $GL_2(\Z_q)$ defined by
  $$
  G_q =
  \AXB{\Z_q\hfill}{\Z_q^*},
  $$
  and set
  $$
  \wt\phi : \axb \eta\xi \in \HP \longmapsto 
  \axb {\phi(\eta)}{\phi(\xi)} \in G_q.
  $$
  Since $G_q$ is a finite group we have that $\Ker(\wt\phi)$ is a
normal subgroup of $\HP$ of finite index.  

  Recall that a while ago we concluded that 
the element $h=\axb\eta\xi$ (introduced near the begining of this proof) 
lies in $\HP\cap\HP^{x\inv}$ if and only if \lcite{\LieInHPxinv}
holds.  We claim that this is the case for all elements
$h\in\Ker(\wt\phi)$.  In fact, if $h\in\Ker(\wt\phi)$, we have that
$\phi(\eta)=0$, and $\phi(\xi)=1$. Therefore there are
$\eta_0,\xi_0\in\AP$, such that
  $
  \eta = q\eta_0,
  $
  and
  $
  \xi = 1+q\xi_0.
  $
  Pluging this in \lcite{\LieInHPxinv} we conclude that
  $$
  \eta a+(1-\xi)b =
  q\eta_0 a-q\xi_0b =
  q\eta_0 {a_0\over q}-q\xi_0{b_0\over q} =
  \eta_0a_0-\xi_0b_0 \in \AP.
  $$  
  This proves that $\Ker(\wt\phi)\subseteq \HP\cap\HP^{x\inv}$, and
hence the index of the latter group in $\HP$ is finite.
  \proofend

Observe that if $\P$ is the empty set then $\AP=\Rat$, and hence
$\HP=G$.  In all other cases we have:

  \state  Proposition 
  \relabel\ExampleNotSubnormal
  If\/ $\P$ is a nonempty set of primes then $\HP$ is not  subnormal in $G$.

  \proof We will show that there exists $h,k\in\HP$, and $x\in G$ such that 
  $$
  (x\inv hx)\inv k (x\inv h x) \notin\HP,
  $$
  thus violating \lcite{\DefSubnormal}.
  Let $a\in\Rat$ and put
  $$
  x=\axb 0a,\quad
  h=\axb 11 \and
  k=\axb 0{-1}.
  $$
  Then
  $$
  x\inv h x = 
  \AXB 0{a\inv} \axb 11 \AXB 0a =
  \axb{a}{1}.
  $$
  So
  $$
  (x\inv hx)\inv k (x\inv h x) =
  \axb{-a}{1} \axb 0{-1} \axb{a}{1} =
  \axb{2a}{-1}.
  $$
  Choosing $a=1/2p$, where $p$ is any prime in $\P$, we conclude that
this is not in $\HP$.
  \proofend

  Among these Hecke pairs we can identify at least one for which $\HP$
is protonormal.

  \state  Theorem
  \relabel\ExampleIsProtonormal
  If\/ $\P=\{2\}$, that is, $\P$ consists of the single prime 2,  then 
  $\HP$ is protonormal in $G$.

  \proof
  Given $x=\axb ba$ in $G$ we need to prove that $\HP^x$ commutes with
$\HP$.  For this let $h=\axb \eta\xi$ and $k=\axb{\nu}{\mu}$ be in
$\HP$ and notice that by \lcite{\ComputeConjugate} we have that
  $$
  x\inv h x k = 
  \AXB {\eta a+(1-\xi )b}{\xi } \AXB{\nu}{\mu} =
  \AXB {\nu+\eta a\mu+(1-\xi )b\mu}{\xi \mu}.
  $$
  We want to to write this as $k'x\inv h'x$, where 
  $h'=\axb {\eta'}{\xi'}$ and $k'=\axb {\nu'}{\mu'}$ are in $\HP$.  We have
  $$
  k'x\inv h'x =
  \AXB {\nu'}{\mu'}  \AXB {\eta'a+(1-\xi')b}{\xi'} =
  \AXB {\eta'a+(1-\xi')b +\nu'\xi'}{\mu'\xi'}.
  $$
  Thus, given $\eta ,\nu\in \AP$, and $\xi ,\mu\in \AP^*$, we need to
find $\eta',\nu'\in \AP$ and $\xi',\mu'\in \AP^*$ such that
  $$
  (\star)\left\{\matrix{\nu+\eta a\mu+(1-\xi )b\mu =  \eta'a+(1-\xi')b +\nu'\xi' \cr\cr
  \xi \mu = \mu'\xi'\hfill}\right.
  $$

  \noindent CLAIM: Seting $\xi'=1+(\xi -1)\mu$, we have that $\xi'\in \AP^*$.

\medskip\noindent
In fact, let $\xi ={x/y}$ and $\mu = z/w$, where $x,y,z$ and $w$ are odd
integers.  Then 
  $$
  \xi' =
  1+(\xi -1)\mu =
  1 + \left({x\over y} - 1\right){z\over w} =
  1 + \left({x-y\over y}\right){z\over w} =
  {yw + (x-y)z\over yw}.
  $$
  Notice that $yw$ is odd and $(x-y)z$ is even so $yw + (x-y)z$ is odd,
hence proving the claim.

  In order to solve $(\star)$ it is then enough to set
  $$
  \def\ccr{\vrule depth 7pt width 0pt\cr}
  \left\{\matrix{
  \xi' =   1+(\xi -1)\mu,\ccr
  \mu' = \xi \mu{\xi'}\inv,\hfill\ccr
  \eta' = \eta \mu,\hfill\ccr
  \nu' = \nu{\xi'}\inv. \hfill
  }\right.
  $$
  \proofend

  \endgroup

  \bigskip \bigskip  \font \sc = cmcsc8 \sc

  \noindent Departamento de Matem\'atica,
  Universidade Federal de Santa Catarina

  \noindent  88040-900 --  Florian\'opolis --
  Brasil 

  \noindent  {exel@\kern1pt mtm.ufsc.br}
  \end